\documentclass[11pt]{article}

\usepackage{amscd,amsmath, amssymb, fancyhdr, mathbbol, dutchcal, url, subfiles}
\usepackage[matrix,arrow,curve]{xy}
\usepackage[backref=page]{hyperref}
\usepackage[dvips]{graphicx}

\renewcommand*{\backrefalt}[4]{%
	\ifcase #1 (Not cited.)%
	\or        (Cited on page~#2.)%
	\else      (Cited on pages~#2.)%
	\fi}

\hypersetup{
	colorlinks   = true,
	citecolor    = magenta,
	linkcolor    =blue,
	urlcolor     =magenta	
}

\numberwithin{equation}{section}

\newcommand{\version}{version 1.0, June 8, 2026}

\def\eqref#1{(\ref{#1})}
\newcommand{\goth}{\mathfrak}
\newcommand{\g}{{\mathfrak g}}

\newcommand{\arrow}{{\:\longrightarrow\:}}
\newcommand{\Z}{{\Bbb Z}}
\def\C{{\Bbb C}}
\newcommand{\R}{{\Bbb R}}
\newcommand{\Q}{{\Bbb Q}}

\def\1{\sqrt{-1}\:}
\newcommand{\restrict}[1]{{\left|_{{\phantom{\cdot}\!\!}_{#1}}\right.}}
\newcommand{\cntrct}                
{\hspace{2pt}\raisebox{1pt}{\text{$\lrcorner$}}\hspace{2pt}}

\renewcommand{\tilde}{\widetilde}
\renewcommand{\bar}{\overline}
\renewcommand{\phi}{\varphi}
\renewcommand{\epsilon}{\varepsilon}
\renewcommand{\geq}{\geqslant}
\renewcommand{\leq}{\leqslant}

\newcommand{\im}{\operatorname{im}}

\newcommand{\Tot}{\operatorname{Tot}}

\newcommand{\diam}{\operatorname{\text{\sf diam}}}

\newcommand{\Area}{\operatorname{\text {\it \sffamily Area}}}

\newcommand{\Hol}{\operatorname{Hol}}

\newcommand{\Lie}{\operatorname{Lie}}

\newcommand{\Diff}{\operatorname{Diff}}

\newcommand{\vce}{{\operatorname{\sf vce}}}

\newcommand{\hor}{{\operatorname{\sf hor}}}

\newcommand{\rk}{\operatorname{rk}}

\newcounter{Mycounter}[section]
\newcounter{lemma}[section]
\renewcommand{\thelemma}{{Lemma \thesection.\arabic{lemma}}}
\newcommand{\lemma}{%
    \setcounter{lemma}{\value{Mycounter}}
    \refstepcounter{lemma}
    \stepcounter{Mycounter}
    {\noindent \bf \thelemma:\ }}

\newcounter{claim}[section]
\renewcommand{\theclaim}{{Claim \thesection.\arabic{claim}}}
\newcommand{\claim}{%
    \setcounter{claim}{\value{Mycounter}}
    \refstepcounter{claim}
    \stepcounter{Mycounter}
    {\noindent \bf \theclaim:\ }}

\newcounter{sublemma}[section]
\renewcommand{\thesublemma}{{Sublemma \thesection.\arabic{sublemma}}}
\newcommand{\sublemma}{%
    \setcounter{sublemma}{\value{Mycounter}}
    \refstepcounter{sublemma}
    \stepcounter{Mycounter}
    {\noindent \bf \thesublemma:\ }}

\newcounter{corollary}[section]
\renewcommand{\thecorollary}{{Corollary \thesection.\arabic{corollary}}}
\newcommand{\corollary}{%
    \setcounter{corollary}{\value{Mycounter}}
    \refstepcounter{corollary}
    \stepcounter{Mycounter}
    {\noindent \bf \thecorollary:\ }}

\newcounter{theorem}[section]
\renewcommand{\thetheorem}{{Theorem \thesection.\arabic{theorem}}}
\newcommand{\theorem}{%
    \setcounter{theorem}{\value{Mycounter}}
    \refstepcounter{theorem}
    \stepcounter{Mycounter}
    {\noindent \bf \thetheorem:\ }}

\newcounter{conjecture}[section]
\newcounter{proposition}[section]
\renewcommand{\theproposition}
      {{Proposition \thesection.\arabic{proposition}}}
\newcommand{\proposition}{%
    \setcounter{proposition}{\value{Mycounter}}
    \refstepcounter{proposition}
    \stepcounter{Mycounter}
    {\noindent \bf \theproposition:\ }}

\newcounter{definition}[section]
\renewcommand{\thedefinition}
      {{Definition~\thesection.\arabic{definition}}}
\newcommand{\definition}{%
    \setcounter{definition}{\value{Mycounter}}
    \refstepcounter{definition}
    \stepcounter{Mycounter}
    {\noindent \bf \thedefinition:\ }}

\newcounter{example}[section]
\renewcommand{\theexample}{{Example \thesection.\arabic{example}}}
\newcommand{\example}{%
    \setcounter{example}{\value{Mycounter}}
    \refstepcounter{example}
    \stepcounter{Mycounter}
    {\noindent \bf \theexample:\ }}

\newcounter{remark}[section]
\renewcommand{\theremark}{{Remark \thesection.\arabic{remark}}}
\newcommand{\remark}{%
    \setcounter{remark}{\value{Mycounter}}
    \refstepcounter{remark}
    \stepcounter{Mycounter}
    {\noindent \bf \theremark:\ }}

\newcounter{problem}[section]
\newcounter{question}[section]
\renewcommand{\thequestion}{{Question \thesection.\arabic{question}}}
\newcommand{\question}{%
    \setcounter{question}{\value{Mycounter}}
    \refstepcounter{question}
    \stepcounter{Mycounter}
    {\noindent \bf \thequestion:\ }}

\newcommand{\proof}{\noindent{\bf Proof:\ }}

\newcommand{\pstep}{\noindent{\bf Proof. Step 1:\ }}

\makeatletter

\@addtoreset{equation}{section} \@addtoreset{footnote}{section}

\def\x@arrow{\DOTSB\Relbar}
\def\xlongrightarrowfill@{\arrowfill@\relbar\relbar\longrightarrow}
\newcommand{\xlongrightarrow}[2][]{%
        \ext@arrow 0099\xlongrightarrowfill@{#1}{#2}}
\makeatother

\def\blacksquare{\hbox{\vrule width 5pt height 5pt depth 0pt}}
\def\endproof{\blacksquare}

\begin{document}

\begin{center}
{\LARGE\bf
Non-commutative Barge-Ghys quasimorphisms: an erratum\\[4mm]
}

Michael Brandenbursky,
Misha Verbitsky

\end{center}

{\small \hspace{0.10\linewidth}
\begin{minipage}[t]{0.85\linewidth}
{\bf Abstract.} This is an erratum to the paper
``Non-commutative Barge-Ghys quasimorphisms''.
The effective version of Ambrose-Singer theorem used
in this paper was mis-stated, making Theorem 3.11 false
as stated. However, this statement remains true 
for a semidirect product $G=A\rtimes K$ of a non-compact
Lie group with a compact Lie group  $K$, as long
as the curvature of the connection belongs to the
Lie algebra of $K$. Theorems 3.12 and 3.22
remain true under this assumptions, the statements
of Theorems 4.28 and 4.31 are modified accordingly.
Theorem 5.2 is valid only for compact Lie groups. 
\end{minipage}
}


\section{Effective Ambrose-Singer theorem}


The proof of Theorem 3.11 in \cite{_BV:quasimorphisms_orig_}
is based on the following observation, called
``effective Ambrose-Singer theorem''. Let $(B, \nabla)$
be a principal bundle with a connection. Then 
its holonomy along a contractible path 
is linearly expressed in terms of the integral 
of the curvature over a disk $D$ filling this path
(\cite[Theorem 1]{_RW:curvature_holonomy_},
\cite{_Deane_Young:holonomy_}).

This observation is true, but we used it to 
obtain a false conclusion: we claimed
that if the curvature has a universal bound, 
then the norm of holonomy is bounded by
a constant times the area of $G$.

This conclusion is false, because, before we 
integrate the curvature, we need to transport
the curvature elements to the origin
using the holonomy of $\nabla$.
Since the path which is used to transport
the curvature elements is not bounded, 
the resulting curvature term is conjugated
by the holonomy action. Unless the adjoint
action of $G$ on the curvature elements is unitary,
this can result in an integral which is unbounded.

To preserve the validity of our paper,
we are forced to consider a Lie group $G=A\rtimes K$, obtained
as a  a semidirect product of a non-compact
Lie group $A$ with a compact Lie group  $K$, and assume
that the curvature of the connection belongs to the
Lie algebra of $K$. In that case, the adjoint
action preserves the norm of the curvature elements,
and the effective Ambrose-Singer theorem can be
applied to obtain the results we announced
in \cite{_BV:quasimorphisms_orig_}.

The last section of \cite{_BV:quasimorphisms_orig_},
generalizing Kazhdan's results 
about $\epsilon$-representations (\cite{_Kazhdan:epsilon_}) 
to an arbitrary Lie group $G$, remains valid as long
as $G$ is assumed compact (or a semidirect product as above).


\section{Corrections to the main results}


In \cite{_BV:quasimorphisms_new_}, we give
corrected versions of the main results of 
\cite{_BV:quasimorphisms_orig_}.

\hfill

{\bf Theorem 4.31, \cite{_BV:quasimorphisms_orig_} (false):}
Let $G$ be a simply connected, connected, non-abelian
rational real nilpotent Lie group, and 
$\Gamma:= \pi_1(M)$,  where $M$ is a closed
manifold of strictly negative sectional curvature, $\dim_\R M >1$.
Then there exists a non-constructible HBG-quasimorphism
$q_\nabla:\; \Gamma \arrow G$.

\hfill

{\bf Theorem 4.14, \cite{_BV:quasimorphisms_new_} (corrected version):}
Let $G := A \rtimes K$ be a semidirect product
of a compact group $K$ and an abelian 
Lie group $A=\R^n, n\geq 2$, and $K$ acts transitively on a sphere in $A$. 
Let $\Gamma:= \pi_1(M)$,  where $M$ is a closed
manifold of strictly negative sectional curvature, $\dim_\R M >1$.
Then there exists a non-constructible HBG-quasimorphism
$q_\nabla:\; \Gamma \arrow G$.

\hfill

Theorem 3.11 and 3.12, \cite{_BV:quasimorphisms_orig_} 
remain valid under an extra assumption:
 $G= A \rtimes K$ is a semidirect product
of a Lie group $A$ with a compact Lie group $K$,
and the curvature 2-form of $\nabla$ takes values
in the Lie algebra $\Lie(G)$ of $K$.

Remark 4.33 of \cite{_BV:quasimorphisms_orig_} is invalid.

\hfill

{\bf Acknowledgements:} We are grateful to the authors of 
\cite{_DFHM:quasimorphisms_to_linear_alg_grps_} for finding
counterexamples to the results of \cite{_BV:quasimorphisms_orig_}.

{\scriptsize

\noindent {\sc Michael Brandenbursky\\
{\sc Ben Gurion University of the Negev\\ 
		Beer Sheva, Israel} \\
{	\tt brandens@bgu.ac.il}\\[4mm]
Misha Verbitsky\\
{\sc Instituto Nacional de Matem\'atica Pura e
	Aplicada (IMPA) \\ Estrada Dona Castorina, 110\\
	Jardim Bot\^anico, CEP 22460-320\\
	Rio de Janeiro, RJ - Brasil }}\\
{	\tt verbit@impa.br}}

\newpage

\ \ 

\vfill

The arXiv.org, in its ineffable 
AI stupidity, decrees that we include the original article. 
We tried to argue, but one does not
argue with the royals (the arxiv admins just ignored us,
sending the same canned replies over and over). 
Apparently arXiv.org is instititionally brain dead now,
and we cannot do anything about it.

The original file is included as is. Please disregard what follows in favor of \url{http://verbit.ru/TeX/New-quasimorphisms.pdf} or \url{https://www.math.bgu.ac.il/~brandens/New-quasimorphisms.pdf}

\vfill

\newpage
\setcounter{tocdepth}{-1}

\begin{center}
{\LARGE\bf
Non-commutative Barge-Ghys quasimorphisms\\[4mm]
}

Michael Brandenbursky,\footnote{M. B. acknowledges support of 
Israel Science Foundation grant 823/23 and 
Humboldt foundation.}
Misha Verbitsky,\footnote{M. V. acknowledges support of 
HSE University basic research program, 
FAPERJ E-26/202.912/2018 and CNPq - Process 310952/2021-2. \\

{\small {\bf 2020 Mathematics Subject
Classification:} 53C05, 20F67, 22E25
}}

\end{center}

{\small \hspace{0.10\linewidth}
\begin{minipage}[t]{0.85\linewidth}
{\bf Abstract.} 
A (non-commutative) Ulam quasimorphism
is a map $q$ from a group $\Gamma$ 
to a topological group $G$ such that
$q(xy)q(y)^{-1}q(x)^{-1}$ belongs
to a fixed compact subset of $G$.
Generalizing the construction 
of Barge and Ghys, we build
a family of quasimorphisms
on a fundamental group of a closed
manifold $M$ of negative sectional curvature,
taking values in an arbitrary Lie group.
This construction, which generalizes
the Barge-Ghys  quasimorphisms, associates a 
quasimorphism to any principal $G$-bundle with 
connection on $M$.

\ \ \ \ \ Kapovich and Fujiwara have shown that
all quasimorphisms taking values in 
a discrete group can be constructed
from group homomorphisms and quasimorphisms 
taking values in a commutative group. 
We construct Barge-Ghys type quasimorphisms 
taking prescribed values on a given subset
in $\Gamma$, producing counterexamples
to the Kapovich and Fujiwara theorem
for quasimorphisms taking values in a Lie group.
Our construction also generalizes a result
proven by D. Kazhdan in his paper 
``On $\epsilon$-representations''.
Kazhdan has proved that for any $\epsilon >0$,
there exists an $\epsilon$-representation of
the fundamental group of a Riemann
surface of genus 2 which cannot be 
$1/10$-approximated by a representation.
We generalize his result by constructing an $\epsilon$-representation
of the fundamental group of a closed manifold
of negative sectional curvature taking values
in an arbitrary Lie group. 
\end{minipage}
}

\tableofcontents


\section{Introduction}


\subsection{The main results}

We start by citing the main results of this paper.
We introduce the relevant definitions and motivation
for these notions later, and here we just give the
statements. The following theorem complements the
main result of \cite{_Fujiwara_Kapovich_} by finding
a counterexample to their theorem when the quasimorphism
takes values in a Lie group.

\hfill

{\bf \ref{_non_constr_Theorem_}:}
Let $G$ be a simply connected, connected, non-abelian
rational real nilpotent Lie group, and 
$\Gamma:= \pi_1(M)$,  where $M$ is a closed
manifold of strictly negative sectional curvature, $\dim_\R M >1$.
Then there exists a non-constructible HBG-quasimorphism
$q_\nabla:\; \Gamma \arrow G$.

\hfill

The HBG-quasimorphisms are defined in 
\ref{_HBG-Definition_} (see also \ref{_G-bundle_HBG_Remark_}),
and constructibility in \ref{_constructible_Definition_}, 
following \cite{_Fujiwara_Kapovich_}. 

\hfill

Another theorem generalizes one of the results
of Kazhdan, \cite{_Kazhdan:epsilon_}.

\hfill

{\bf \ref{_no_approx_epsilon_rep_Theorem_}:}
Let $M$ be a closed manifold
of strictly negative sectional curvature, $G$ a positive-dimensional
connected Lie group, and 
$P$ a trivial principal $G$-bundle. For any connection $\nabla$
in $P$, let $q_\nabla:\; \pi_1(M) \arrow G$ denote
the corresponding HBG-quasimorphism
(\ref{_HBG-Definition_}). Choose a
left-invariant metric on $G$ such that
the diameter of any closed subgroup is at least $1/3$.
Then for each $\epsilon>0$,
there exists a connection $\nabla$ such that
$q_\nabla$ is an $\epsilon$-representation
which cannot be $1/3$-approximated by 
a representation.

\hfill

We explain the terms
``$\epsilon$-representation''
and ``$\delta$-approximation''
in Subsection 
\ref{_approximation_Subsection_}, following \cite{_Kazhdan:epsilon_}.

\subsection{Ulam quasimorphisms}

Throughout this paper, $G$ means a Lie group
of algebraic type; usually we assume that $G$
is connected, though this assumption is not
always necessary.

\hfill

\definition\label{_Ulam_qm_Definition_}
Let $G$ be a Lie group, and $\Gamma$ any group. 
 {\bf An Ulam quasimorphism}
$q:\; \Gamma \arrow G$ is a map which satisfies
$q(x^{-1}) = q(x)^{-1}$ and  $q(x) q(y) \in K\cdot q(xy)$,
where $K\subset G$ is a fixed compact,
independent from the choice of $x, y$.
Two Ulam quasimorphisms $q_1, q_2:\;\Gamma \arrow G$ are called
{\bf equivalent} if there exists a compact subset $K\subset G$ 
such that $q_1(x) \in K\cdot q_2(x)$, for all $x\in \Gamma$.

\hfill

In this paper, we give a versatile construction of 
Ulam quasimorphisms associated with connections
in vector bundles over a closed manifold $M$ of 
strictly negative sectional curvature. As a result, we
obtain a quasimorphism from $\pi_1(M)$ to a Lie group,
called {\em a Barge-Ghys quasimorphism.}
This allows us to find quasimorphisms which are
not constructible, in the sense of Fujiwara and Kapovich
\cite{_Fujiwara_Kapovich_}. The 
Fujiwara-Kapovich constructible quasimorphisms
are obtained from quasimorphisms taking
values in abelian groups. Unlike the constructible quasimorphisms,
the Barge-Ghys quasimorphisms we construct cannot be  
obtained from ones taking values in abelian groups (\ref{_non_constr_Theorem_}).

Interestingly enough, Fujiwara-Kapovich 
have proved that {\em any} Ulam quasimorphism taking
values in a discrete group is constructible.
It turns out that the Ulam quasimorphisms taking
values in Lie groups are of entirely different nature.

Ulam quasimorphisms from a free group to a Lie group
$G$ with bi-invariant metric were explored
by P. Rolli in \cite{_Rolli:free_quasimorphisms_}.
Rolli has constructed non-trivial quasimorphisms from a free
group to $G$. It seems that Rolli's construction also
can lead to non-constructible quasimorphisms.

As an application of our approach, we give a  
construction of $\epsilon$-repre\-senta\-tions which cannot
be approximated by a representation 
(\ref{_no_approx_epsilon_rep_Theorem_}), generalizing
a result of Kazhdan (\cite{_Kazhdan:epsilon_}).

Note that our definition of Ulam quasimorphism
may not be the optimal for some purposes. Fujiwara and Kapovich 
\cite{_Fujiwara_Kapovich_} give several non-equivalent
(more relaxed) definitions of a quasimorphism taking values in 
a non-abelian group $G$; for abelian $G$, all these
definitions are equivalent. Another, even more relaxed, definition
was considered in \cite{_Hartnick_Schweitzer_}.
We give a more detailed presentation of these
notions in Subsection \ref{_alge_geome_ulam_Subsection_}.

One of the results we obtain builds on the difference
between these notions. Recall that the {\bf geometric
quasimorphism} (\cite{_Fujiwara_Kapovich_})
is a map $q:\; \Gamma \arrow G$ such that 
there exists a compact subset $K\subset G$ 
such that $q(xy) \in K q(x) K q(y)$ for all $x, y\in \Gamma$. 

\hfill

\theorem
Let $\Gamma$ be a fundamental group of 
a closed manifold of strictly negative curvature,
$G$ a non-commutative, simply connected
nilpotent Lie group, and $\Lambda$ a cocompact
lattice in $G$.\footnote{By Maltsev's theorem 
(\cite{_Corwin_Greenleaf_}), existence of a
cocompact lattice is equivalent to $G$ being rational.}
Then a non-constructible
HBG-quasimorphism $q:\; \Gamma \arrow G$
obtained in \ref{_non_constr_Theorem_}
can be approximated by a geometric quasimorphism
$q_0:\; \Gamma \arrow \Lambda$, which is also
non-constructible.

\proof \ref{_geometric_not_Ulam_Remark_}. \endproof

\hfill

Notice that, by contrast, any 
{\em Ulam} quasimorphism $q_1:\; \Gamma \arrow \Lambda$
is by \cite{_Fujiwara_Kapovich_} constructible.
This is where the difference between the Ulam
quasimorphisms and the geometric quasimorphisms
becomes apparent.


\subsection{Quasimorphisms, bounded cohomology and the commutator length}

In the literature, the quasimorphisms are
usually considered as maps taking values in $\R$.
In this context, a quasimorphism is 
defined as a map $q:\; G \arrow \R$ 
which satisfies $|q(xy) -q(x)-q(y)|< C$,
where $C$ is a constant independent from $x, y$. 
We shall sometimes call such maps {\bf commutative
quasimorphisms}.

In geometruc context, this notion originates in the paper of Gromov
\cite{_Gromov:Bounded_}. Using the quasimorphisms, R. Brooks
proved that the second bounded
cohomology $H^2_b({\Bbb F}_2,\R)$ of the free group ${\Bbb F}_2$
is infinite-dimensional (\cite{_Brooks:bounded_,_Gromov:Bounded_}).
Since then, quasimorphisms $q:\; G \arrow \R$
became prevalent in topology (\cite{_Calegari:scl_}),
symplectic geometry (\cite{_EP:Calabi_qm_,_Shelukhin:thesis_})
and dynamics (\cite{_Brandenbursky_Marcinkowski_,_Gambaudo_Ghys_}).

Barge and Ghys (\cite{_Barge_Ghys:bounded_}) generalized the
observation of Brooks to prove that the second bounded cohomology 
of $\pi_1(S)$ is infinite-dimensional
for any Riemann surface $S$ with $g(S) >1$.

We say that two quasimorphisms $q, q_1:\; G \arrow \R$ are {\bf equivalent}
if $|q(x) - q_1(x)| < C$, where $C$ is a constant independent from $x$.
The group ${\cal Q}$ of equivalence classes of quasimorphisms
fits into an exact sequence
\begin{equation}\label{_bounded_exact_Equation_}
H^1(G, \R) \arrow {\cal Q} \arrow H^2_b(G, \R) \arrow H^2(G,\R),
\end{equation}
where $H^2_b(G, \R)$ denotes the bounded cohomology.
A quasimorphism $q:\; G \arrow \R$ is called {\bf homogeneous}
if  it satisfies $q(x^n) = n q(x)$ for any $x\in G$ and $n\in \Z$.
It is possible to see that
every quasimorphism $q:\; G \arrow \R$  is equivalent to a unique
homogeneous quasimorphism (\cite[Lemma 2.2.1]{_Calegari:scl_}).

The commutator length of $g\in [G,G]$ is the
minimal number $m$ such that $g$ can be
represented as a product of $m$ commutators.
Recall that {\bf the stable commutator length}
of an element $g\in G$ is defined as
\[
\mathop{\sf scl}(g):= \lim_{n\to\infty} \frac{\mathop{\sf cl}(g^n)}n,
\]
where $\mathop{\sf cl}$ is the commutator length. 
The homogeneous quasimorphisms are related with the
stable commutator length, due to the celebrated theorem
of Bavard (\cite{_Bavard:scl_}). This result is known
as {\em Bavard duality}. In its simplest version,
Bavard duality can be stated as follows: 
for any $g\in [G,G]$,
$g$ has non-zero stable commutator length if and only if
there exists a homogeneous quasimorphism  $q:\; G \arrow \R$
such that $q(g)\neq 0$. 


\subsection{Barge-Ghys construction and manifolds of strictly negative curvature}

Further on, we are interested in vector bundles over
closed manifolds of strictly negative sectional curvature.
Sometimes such manifolds are called ``hyperbolic'', but
we don't use this term to avoid confusion with manifolds
of constant strictly negative curvature.

We study the Ulam quasimorphisms associated with the
holonomy of a connection on vector bundles or on principal
$G$-bundles, where $G$ is a connected Lie group. In the sequel,
we use either the language of $G$-bundles or the language
of vector bundles with connection, whatever is more convenient.
However, all statements that we make can be easily translated
from one language to another, giving equivalent results in
two parallel frameworks.

\hfill

\remark\label{_dim>1_Remark_}
We tacitly assume that the base manifold
of strictly negative curvature has dimension at least 2. 

\hfill

Barge-Ghys quasimorphisms were introduced in
\cite{_Barge_Ghys:bounded_}, who used them to prove
that the bounded cohomology of a Riemann
surface is infinite-dimensional. In the literature
these quasimorphisms are variously called de Rham quasimorphisms
(\cite{_Calegari:scl_}) or Barge-Ghys quasimorphisms
\cite{_EP:Calabi_qm_,_RP:functional_}; we follow
the second convention.

The generalization of Barge-Ghys quasimorphisms
to closed manifolds of strictly negative curvature
seems to be well known (\cite{_Marasco:Cup_}).
However, the infinite-dimensionality
of the space of Barge-Ghys quasimorphisms,
originally proven in \cite{_Barge_Ghys:bounded_} 
for Riemann surface, is less straightforward.
This result was established in \cite{_BFMSS:bounded_}.

\subsection{Barge-Ghys quasimorphisms associated with a connection}

We feel compelled to make a few 
comments about the terminology used in this paper.

Originally, the Barge-Ghys quasimorphisms were
associated with a non-closed 1-form $\theta$ on a closed manifold $M$
of strictly negative sectional curvature.
Given $\gamma\in \pi_1(M, x)$, we represent $\gamma$ by
a geodesic loop $\underline\gamma$ (which is unique, because $M$ has negative
curvature), and put $q(\gamma) = \int_{\underline \gamma} \theta$.

Taking $\theta$ as a connection form in a trivial 
line bundle, the number $\int_{\underline \gamma} \theta$
is intepreted as the holonomy of this connection
along $\underline \gamma$ (Subsection
\ref{_Barge_Ghys_noncomm_Subsection_}).

In Subsection
\ref{_Barge_Ghys_noncomm_Subsection_}, we generalize this
construction to arbitrary vector bundles
and to principal $G$-bundles  with connection.
This gives a ``non-commutative quasimorphism''
(\ref{_Ulam_qm_Definition_}), which we call
``a non-com\-mu\-ta\-tive Barge-Ghys quasimorphism 
associated with a connection''.

This quasimorphism does not satisfy
$q(x^n)= q(x)^n$, that is, it is not homogeneous
(\ref{_homo_Definition_}).  In the usual
theory of quasimorphisms (which we call
throughout this paper ``commutative quasimorphisms'')
for every quasimorphism there exists a unique
homogeneous quasimorphism in the same equivalence
class. This construction is called ``homogenization'';
it is obtained by the standard limit construction 
(\cite[Lemma 2.2.1]{_Calegari:scl_}). 

Unfortunately, we were unable to generalize the homogenization
construction to the non-commutative case. Instead we use
an ad hoc construction, which works in the same
generality, and gives a homogeneous quasimorphism
(Subsection \ref{_homo_BG_Subsection_}).
The quasimorphisms obtained this way are called
``the homogeneous Barge-Ghys quasimorphisms''.
They are also associated with a connection in
a vector bundle or in a principal $G$-bundle.

There are disclaimers we need to make at this point.
First of all, we could not devise a general
definition of ``Barge-Ghys quasimorphisms''.
These two kinds of Barge-Ghys quasimorphisms
are the only cases we could came up with.

Second, the Barge-Ghys quasimorphisms associated with
a connection are distinct from the homogeneous Barge-Ghys quasimorphisms.
These are two distinct classes which rarely intersect, and
we do not know if they can be united in a meaningful
class of more general quasimorphisms. 
However, each Barge-Ghys quasimorphism
accociated with a $G$-bundle with connection
is equivalent to the homogeneous Barge-Ghys quasimorphism
associated with the same connection (\ref{_HBG_equiv_to_BG_Claim_}).

\subsection{Ulam stability and Kazhdan's $\epsilon$-representations}
\label{_approximation_Subsection_}

The earliest mention of quasimorphisms 
is found in Ulam's 1960 book ``A collection of mathematical problems'',
Chapter 6 (\cite{_Ulam:stability_}). Ulam 
defined an $\epsilon$-automorphism of a
topological group as a map $\rho:\; G \arrow G$
such that $\rho(xy) \rho(x)^{-1} \rho(y)^{-1}$
belongs to an $\epsilon$-neighbourhood of the identity.
Ulam asked whether any such map admits
an $k\epsilon$-approximation by an automorphism,
for some $k>0$ which is independent from $\epsilon$.

More recently, this problem was generalized
to representations (see \ref{_Ulam_stability_Definition_} below).
However, the notion of Ulam stability seems to 
be present, implicitly, in earlier works
of von Neumann (\cite{_von_Neumann:Mechanik_}) 
and Turing (\cite{_Turing:Lie_}).

A similar question was considered in 1982 by
D. Kazhdan (\cite{_Kazhdan:epsilon_}). 
He defined {\bf an $\epsilon$-representation}
of a group $G$ as a map $\rho:\; G \arrow U(V)$,
where $V$ is a Hilbert space, finite dimensional or infinite dimensional,
satisfying \[ \|\rho(xy) - \rho(x) \rho(y)\| <\epsilon,\]
where $\|\cdot \|$ is the operator norm.
The distance between two maps 
$\rho_1, \rho_2:\; G \arrow U(V)$ can be
defined as \[ d(\rho_1, \rho_2):= \sup_{x\in G}\|\rho_1(x) -\rho_2(x)\|.\]
Following a suggestion of V. Milman, Kazhdan
asked whether for any $\delta>0$ there exists $\epsilon>0$
such that any  $\epsilon$-representation
can be $\delta$-approximated by a representation.
He proved that this holds true for amenable groups.
When $G=\pi_1(S)$, where $S$ is a genus 2 Riemann surface,
Kazhdan has constructed an $\epsilon$-representation $G \arrow U(n)$, 
which cannot be $1/10$-aproximated by a representation,
for any given $\epsilon >0$. 



\hfill

\definition\label{_Ulam_stability_Definition_}
The group $\Gamma$ is called {\bf Ulam stable}
(\cite{_BOT:Ulam_stability_}) if for any $\delta>0$ there exists
$\epsilon>0$ such that any
finite-dimensional $\epsilon$-representation
 $q:\; \Gamma \arrow U(V)$ 
can be $\delta$-approximated by a representation
$\rho:\; \Gamma \arrow U(V)$. 
It is called
{\bf strongly Ulam stable} if
the same is true even for infinite-dimensional
Hilbert representations.

\hfill

In \cite{_Kazhdan:epsilon_}, D. Kazhdan has
proven that all amenable groups are strongly Ulam stable.
Using the Barge-Ghys quasimorphisms, we were
able to prove this result for the fundamental
group of any closed strictly negatively curved manifold.
In \ref{_no_approx_epsilon_rep_Theorem_}, we generalize 
\cite[Theorem 2]{_Kazhdan:epsilon_}
showing that a fundamental group $\Gamma$
of a complete Riemannian manifold with 
uniformly bounded strictly negative
sectional curvature admits
an $\epsilon$-representation
taking values in any given Lie group $G$, 
which cannot be $1/3$-approximated by
a representation\footnote{Our norm conventions are different
from Kazhdan's, but the constants $1/3$ 
and $1/10$ depend on the choice of normalizations, and 
are not important in themselves.}.
Kazhdan proves this for $\Gamma=\pi_1(S)$,
where $S$ is a genus 2 Riemann surface, and $G=U(n)$.

In \cite{_BOT:Ulam_stability_}, Burger, Ozawa and Thom
address the question of strong Ulam stability, obtaining 
definite results for a large class of groups. 
Let $\Gamma$ be a group which contains
a subgroup $\Lambda$ admitting a homogeneous $\R$-valued quasimorphism 
which is not a homomorphism, that is, such that
the map \[ H^2_b(\Lambda, \R) \arrow H^2(\Lambda,\R)\]
has non-trivial kernel. Then $\Gamma$ has
an infinite-dimensional Hilbert $\epsilon$-repre\-sen\-tation,
for any given $\epsilon >0$, which 
cannot be $\frac{\sqrt 3}{16}$-approximated by
a representation. Also, they observe that
any free group admits a finite-dimensional
$\epsilon$-representation taking values in $U(n)$
which cannot be 2-approximated by a representation 
(\cite[Proposition 3.3]{_BOT:Ulam_stability_}).

In \cite{_GLMR:asymptotic_cohomology_},
the Ulam stability was cast in a cohomological setting.
Glebsky, Lubotzky, Monod and Rangarajan
defined an asymptotic version of bounded
cohomology, and proved that vanishing of
asymptotic cohomology implies Ulam stability.
In the same paper, Ulam stability
of $U(1)$-representations is directly related
to existence of quasimorphisms.

In this paper, we further generalize Kazhdan's theorem,
which is known for $\epsilon$-representations
with values in $U(n)$ (\cite[Proposition 3.3]{_BOT:Ulam_stability_})
to  $\epsilon$-representations
taking values in an arbitrary Lie group
(\cite{_no_approx_epsilon_rep_Theorem_}).

\subsection{Constructible quasimorphisms}

The paper \cite{_Fujiwara_Kapovich_} by
Fujiwara and Kapovich is the fundamental
treatment of Ulam quasimorphisms taking values
in a non-commutative group $G$. Fujiwara and Kapovich 
considered the case when $G$ is discrete, and their
result is mostly negative. Fujiwara and Kapovich 
defined ``constructible quasimorphisms'', which
are up to equivalence and finite index
quasimorphisms which can be constructed in
terms of  homomorphisms, quasimorphisms with abelian
target, and sections of bounded central extensions.

Then they proved
that any Ulam quasimorphism taking values in 
a discrete group is always constructible 
(\cite[Theorem 1.2]{_Fujiwara_Kapovich_}).

In \ref{_non_constr_Theorem_} we  
construct examples of quasimorphisms
taking values in a Lie group $G$, which do not
satisfy conclusions of
\cite[Theorem 1.2]{_Fujiwara_Kapovich_}.
There is no contradiction because the group $G$ is not discrete.
We prove that for any closed manifold $M$ of 
strictly negative sectional curvature, 
there exists an Ulam quasimorphism
$\pi_1(M) \arrow G$ which is not 
constructible.

\section{Principal bundles and vector bundles}
\label{_principal_Section_}

Almost everything we say in this paper can be
formulated for vector bundles or for principal $G$-bundles,
where $G$ is a Lie group. Usually we state only one
of two versions, leaving the rest for the reader.
In this section, we briefly explain the passage
from vector bundles with connection to $G$-bundles 
and vice versa.

\hfill

\definition
Let $G$ be a Lie group.
A {\bf principal $G$-bundle}
is a smooth manifold $E$ equipped with
a smooth, free $G$-action, such that the
natural map $E\to E/G$ is a locally 
trivial fibration.

\hfill

\example
Consider the standard action of $U(1)$ on 
3-dimensional sphere $S^3\subset \C^2$,
with $e^{i t}$ taking $(\xi_1, \xi_2)$
to $(e^{i t}\xi_1, e^{i t}\xi_2)$.
Then $S^3/U(1)= S^2$, and this action
defines a principal $U(1)$-bundle on
a 2-dimensional sphere. This fibration
is known as {\bf Hopf fibration}.

\hfill

\definition
Let $\pi:\; E\arrow M$ be a principal $G$-bundle
and $V$ a space with $G$-action.
Consider the quotient $(E\times V)/G$,
where $G$ acts diagonally.
Since the action of $G$ is free on fibers
of $\pi$, the quotient $(E\times V)/G$
is a locally trivial fibration on $M$
with fiber $V$. It is called {\bf the associated
fibration}.

\hfill

\example
Let $V$ be a representation of a group $G$, and
$\pi:\; E\arrow M$ be a principal $G$-bundle.
Then the associated fibration $(E\times V)/G$
is a vector bundle, called {\bf a vector bundle with
a $G$-structure}. We say that {\bf the structure group
of a vector bundle $B$ is reduced to $G$}
if $B$ is obtained from a principal $G$-bundle this way.

\hfill

\example
Consider a complex vector bundle $B$ over $M$ equipped with
a Hermitian structure, and let $\pi:\; E\arrow M$
be the space of all orthonormal complex frames in $B$.
Since the group $U(n)$ acts on orthonormal complex frames
freely and transitively, $E$ is a principal $U(n)$-bundle. 
Consider the standard
$\C^n$-representation $V$ of $U(n)$. By construction,
the vector bundle $(E\times V)/U(n)$ coincides with $B$,
hence this construction gives a reduction of the structure 
group of $B$ to $U(n)$.

\hfill

We have described the correspondence between the
vector bundles (with appropriate reduction of the structure group)
and the principal $G$-bundles. It turns out that this
construction is well compatible with connections;
one can define the holonomy and the curvature in both
contexts, and these notions are equivalent.
This is a part of standard differential geometry
course, see e. g. \cite{_Sternberg:lectures_,_Kob_Nomizu_}.

We will presently give a partial description of this correspondence.

\hfill

\definition\label{_Ehresmann_conne_Definition_}
Let $\pi:\;E \arrow M$ be a smooth fibration, with
$T_\pi E$ {\bf the bundle of vertical tangent vectors}
(vectors tangent to the fibers of $\pi$).
An {\bf Ehresmann connection} on $\pi$
is a sub-bundle $T_\hor E\subset TE$ such that
$TE= T_\hor E\oplus T_\pi E$.
The {\bf parallel transport} along the path $\gamma:\; [0, a]\arrow Z$
associated with the Ehresmann connection
is a diffeomorphism $V_t:\;  \pi^{-1}(\gamma(0))\arrow  \pi^{-1}(\gamma(t))$
smoothly depending on $t\in [0, a]$ and 
satisfying $\frac {dV_t}{dt}\in T_\hor E$.

\hfill

\definition\label{_linear_Ehresmann_Definition_}
Let $B$ be a vector bundle on $M$ and $\Tot B \stackrel \pi \arrow M$
its total space. An Ehresmann connection
on $\pi$ is called {\bf linear} if it is preserved by 
the homothety map $\Tot B \arrow \Tot B$ 
mapping $v$ to $\lambda v$ and by the addition map
$\Tot (B\oplus B) \arrow \Tot B$,
that is, addition preserves horizontal vectors.

\hfill

\proposition\label{_linear_Ehresmann_Proposition_}
A notion of a linear Ehresmann connection
on a vector bundle $B$ coincide with
the usual notion of a connection;
the corresponding parallel transport
maps coincide as well.

\proof \cite[Proposition 2.21]{_Ornea_Verbitsky:Book_}.
\endproof

\hfill

\definition
Let now $\pi:\;E \arrow M$
be a principal $G$-bundle. A $G$-connection on $\pi$
is a $G$-invariant Ehresmann connection.

\hfill

\remark\label{_product_conne_Remark_}
Let $\pi:\;E \arrow M$ be a principal $G$-bundle
with $G$-connection $\nabla$, and $X$ a smooth manifold with $G$-action $\rho$. 
Then the associated bundle $E_X:=(E\times X)/G$ inherits an Ehresmann connection:
$TE_X = T_\hor E \oplus (TX/\rho(\g))$, where $\g=\Lie(G)$ is the Lie algebra of $G$.

\hfill

\definition
 Let $B$ be a vector bundle on $M$ with structure
group $G$,  $V$ representation of $G$, 
and $E$ the corresponding principal $G$-bundle,
$B=(E\times V)/G$. Then for any  $G$-connection on $E$
induces an Ehresmann connection on $B$ as in \ref{_product_conne_Remark_}.
Clearly, this connection is linear, hence induces
a connection on $B$ as on a vector bundle
(\ref{_linear_Ehresmann_Proposition_}).
This connection is called {\bf induced by the
$G$-connection $\nabla^E$.}

\hfill

The holonomy of a $G$-connection and its induced
connection are compatible, in the following sense:

\hfill

\claim Let $E$ the a principal $G$-bundle on $M$, $V$ a representation of $G$,
and $B=(E\times V)/G$ the corresponding  vector bundle with the structure group $G$.
Let $\nabla^E$ a $G$-connection on $E$, and $\nabla$ the
induced connection on $B$.
 Denote the holonomy of $\nabla^E$ along a loop $\gamma$ 
based on $m\in M$ by $g_\gamma\in G$. Then the holonomy of $\nabla$ along
$\gamma$ is obtained from the action of $g_\gamma$ on 
$B\restrict m = \left(E\restrict m \times V\right)/G$.

\proof \cite{_Kob_Nomizu_}.
\endproof

\hfill

Throughout this paper, we use one of these equivalent languages,
and assume tacitly an analogous statement for the other one.
It is slightly more convenient to speak of non-commutatibe 
Barge-Ghys quasimorphisms in the language of vector bundles:
this way, the correspondence with the usual, commutative Barge-Ghys
quasimorphism is more apparent. However, it is more natural to state 
and prove the generalization of Kazhdan's theorem in the language
of principal $G$-bundles. The translation from one of these
languages to another is straightforward and is left to the reader 
as an exercise.

\section[Barge-Ghys quasimorphisms on fundamental groups
  of closed manifolds of strictly negative curvature]{Barge-Ghys quasimorphisms on fundamental \\groups
  of closed manifolds of strictly negative curvature}

\subsection{Manifolds with strictly negative sectional curvature}

In this preliminary section we list
a few arguments of Riemannian geometry, 
most of them either classical or due to M. Gromov, \cite{_Gromov:Bounded_}.
For manifolds of constant negative curvature,
all the results we are going to obtain are classical and well known
(\cite[Chapter 8]{_Frigerio:bounded_}); however, for arbitrary
Riemannian manifolds of strictly negative curvature, a more subtle
approach is required.

\hfill

\theorem\label{_Cartan_Hadamard_Theorem_}
(Cartan-Hadamard)\\
Let $M$ be a complete, simply connected Riemannian manifold
of non-positive sectional curvature. Then
$M$ is contractible.

\hfill

\proof We give a sketch of the proof,
following \cite{_BBI_}. In \cite{_BBI_}
this result was stated and proven for
CAT(0)-spaces, but the comparison 
inequalities which are required by the
CAT(0)-geometry easily follow from 
 the non-positivity of the sectional curvature.

By Hopf-Rinow theorem, every two
points of $M$ can be connected by a geodesic.
Let $\gamma_1:\;[a,b]\arrow M $ and
$\gamma_2:\;[c,d]\arrow M$  be segments of
geodesics in $M$, parametrized by the arc length.
As shown in \cite{_BBI_}, the distance function
$D:\; [a,b]\times [c,d]\arrow \R^{>0}$
taking $x, y$ to $d(\gamma_1(x), \gamma_2(y))$
is strictly convex, unless the geodesics
$\gamma_1, \gamma_2$ are segments of the same
geodesic line. This implies, in particular,
that any two points are connected by a
unique geodesic: indeed, if $\gamma_1$ and
$\gamma_2$ have the same ends,  $\gamma_1(a)=\gamma_2(c)$
and $\gamma_1(b)=\gamma_2(d)$ the 
function $D$ would be equal to 0
in $(a,c)$ and $(b,d)$, hence
it is zero on the diagonal, and the
images of $\gamma_1$ and $\gamma_2$ coincide.

Fix a reference point $p\in M$ and consider 
the function \[ H:\; M\times[0,1]\arrow M\]
taking $x\in M$ and $t\in[0, 1]$ to 
$\gamma(t\cdot d(p,x))$, where $\gamma:\; [0, d(p,x)]$ 
is the geodesic  connecting $p$ to $x$.
A similar argument implies that
$H$ is continuous; clearly, $H$ is
a deformation retraction of $M$ to $p$,
hence $M$ is contractible.
\endproof

\hfill

We will not use the Cartan-Hadamard theorem,
but we use its corollary, which is inherent
in its proof.

\hfill

\corollary\label{_geode_homoto_Corollary_}
Let $M$ be a simply connected, complete
manifold of non-positive sectional curvature.
Fix a reference point $p\in M$ and consider 
the function $H_p:\; M\times[0,1]\arrow M$
taking $x\in M$ and $t\in[0, 1]$ to 
$\gamma(t\cdot d(p,x))$, where $\gamma:\; [0,
  d(p,x)]\arrow M$ 
is the geodesic  connecting $p$ to $x$.
Then $H_p$ is continuous. \endproof

\hfill

We could use the map $H_p$ to define {\bf a straight singular
simplex,} following \cite{_Gromov:Bounded_}.
Since we need only 2-dimensional simplices,
we restrict ourselves to the 2-dimensional case.

\hfill

\definition
Let $M$ be a simply connected, complete
manifold of non-positive sectional curvature.
Let $\gamma$ be the geodesic segment
connecting $b$ to $c$. We obtain a triangle
by connecting $a$ to all points of $\gamma$ by
a unique geodesic.
The {\bf geodesic simplex} $\Delta(a,b,c)$,
associated with the points $a, b, c\in M$
is the union $\cup_{t\in [0,1]} [H_a(t)(b),H_a(t)(c)]$,
where $H_a(t):\; M \arrow M$ is
the homotopy along geodesics passing through $a$, defined in 
\ref{_geode_homoto_Corollary_}, and
$[H_a(t)(b),H_a(t)(c)]$ the geodesic segment connecting
$H_a(t)(b)$ and $H_a(t)(c)$.

\hfill

\remark
The boundary of the simplex $\Delta(a,b,c)$
is the union of geodesics connecting
$a$ to $b$, $b$ to $c$ and $c$ to $a$.
Indeed, the homothety $H_a(t)$
moves any point $x\in M$ along the
geodesic connecting $x$ to $a$
as $t$ goes from $1$ to $0$,
hence the segments of the boundary
connecting $b$ to $a$ and $c$ to $a$
are geodesics; the third segment
is $\gamma$, which is also chosen
geodesic.

\hfill

\remark
Note that the order of the points
$a, b, c\in M$ is important. Indeed,
unless $M$ has constant sectional curvature,
the simplex $\Delta(a,b,c)$ and 
(say) $\Delta(b,a,c)$ are different:
otherwise, if $\Delta(a,b,c)=\Delta(b,a,c)=\Delta(a,c,b)$,
this simplex is a segment of a completely geodesic
2-dimensional submanifod, and a general 
Riemannian manifold does not have completely geodesic
submanifolds (\cite{_Murphy_Wolhelm:no_geode_subma_}). 

\hfill

The main technical result about manifolds
of strictly negative sectional curvature which we use
is the following theorem of Gromov, which is
proven in \cite{_Gromov:Bounded_}
for any straight singular simplex.

\hfill

\theorem\label{_Gromov_volume_triangle_Theorem_}
Let $M$ be a simply connected, complete
manifold of strictly negative sectional curvature 
$K(M) <-\epsilon <0$, and $\Delta(a,b,c)$ a geodesic
simplex defined above. Then 
the Riemannian area $\Area(\Delta(a,b,c))$ 
satisfies  $\Area(\Delta(a,b,c)) \leq \pi\epsilon^{-2}$.

\proof \cite[\S I.3]{_Gromov:Bounded_}. \endproof

\hfill

In the sequel, we will need two statements about uniqueness of geodesics
on manifolds of strictly negative sectional curvature. The first is a
direct consequence of the Cartan-Hadamard's theorem (\ref{_Cartan_Hadamard_Theorem_}).

\hfill

\claim
Let $\gamma\in \pi_1(M,p)$ be an element
of a fundamental group of a closed  manifold of strictly negative sectional curvature.
Then $\gamma$ is represented by a unique geodesic loop
based at $p$.
\endproof

\hfill

Another statement (slightly less trivial) deals with
{\em free geodesic loops}. Recall that a {\bf free geodesic
loop} is an immersed submanifold of dimension 1 which is locally geodesic.

\hfill

\proposition\label{_free_geodesic_unique_Proposition_}
Let $M$ be a closed  manifold of strictly negative sectional curvature,
and $\phi:\; S^1 \arrow M$ be a smooth map.
Then there exists a unique free geodesic loop $\phi_1$
which is free homotopic to $\phi$. Moreover, $\phi_1$
strictly minimizes the length of the loop.

\proof \cite[Theorem 3.8.14]{_Klingenberg:Riemannian_}. 
\endproof

\subsection{Non-commutative Barge-Ghys quasimorphisms}
\label{_Barge_Ghys_noncomm_Subsection_}

Throughout this section, $G$ is a connected Lie group.

\hfill

\definition\label{_NC_BG_Definition_}
Let $M$ be a closed manifold with non-positive sectional curvature,
and $(P, \nabla)$ a principal $G$-bundle
 with connection. Fix $x\in M$. The
{\bf non-commutative Barge-Ghys map}
takes $\gamma \in \pi_1(M, x)$  
to the holonomy of $\nabla$ along
the geodesic path starting and ending at $x$
and homotopic to $\gamma$.\footnote{Since $M$ has
  non-positive curvature, this geodesic path is unique
in its homotopy class, \ref{_free_geodesic_unique_Proposition_}.}

\hfill

\remark 
Throughout this paper, we consider the maps
$q:\; \Gamma \arrow G$, where $\Gamma=\pi_1(M)$ is the
fundamental group of a closed manifold of strictly 
negative curvature. Such groups are Gromov hyperbolic,
however, it is not clear yet whether our constructions
can be generalized to all Gromov hyperbolic groups.
Nevertheless, all results we prove remain valid
whenever each geodesic simplex in the universal
cover of $M$ can be filled by a disk of bounded 
area, as in \ref{_Gromov_volume_triangle_Theorem_}.
Note that the homological version of the
bounded filling result is also true for hyperbolic
groups which are not of geometric origin
(\cite{_Mineyev:bounded_}).

\hfill

\theorem\label{_geodesic_polygon_holonomy_Theorem_}
Let $M$ be a closed manifold of strictly negative
sectional curvature, and $\Theta$ a 
geodesic $n$-polygon in $M$, that is, a contractible
loop of $n$ geodesic segments. Consider
a principal $G$-bundle $(P, \nabla)$ with connection on
$M$, and let $h(\Theta)\in G$ 
be the holonomy along the boundary of $\Theta$,
considered as a loop starting and ending at $p\in \Theta$.
Then $h(\Theta)$ belongs to a compact $K_n\subset G$
which is the same for all $n$-polygons $\Theta$,
but depends on $n$ and $(P, \nabla)$
and the bound on the curvature of $M$.

\hfill

\proof
By the effective version of the Ambrose-Singer theorem,
the holonomy along a path is linearly expressed in terms of the integral 
of the curvature over a disk filling this path
(\cite[Theorem 1]{_RW:curvature_holonomy_},
\cite{_MO:holonomy_question_},
\cite{_Deane_Young:holonomy_}).
Therefore, for a left-invariant metric
on $G$, the holonomy is bounded in terms of the
integral of the curvature.
The area of any geodesic simplex is
bounded by \ref{_Gromov_volume_triangle_Theorem_}.
The absolute value of the
curvature of $\nabla$ is bounded from above because
$M$ is compact, and the curvature form on the pullback bundle 
$(\tilde P, \tilde \nabla)$
is obtained as a pullback of the
curvature of $(P, \nabla)$.

This implies that $h(\Theta)$ belongs to a fixed compact 
$K$ when $n=3$ and $\Theta$ is a simplex. When $n >3$, 
we represent $\Theta$ as a boundary of the union
of $n-2$ geodesic simplexes $D_1, .., D_{n-2}$
with common vertex $p$. Then the holonomy $h(\Theta)$
is obtained as a product $h(\Theta)=h(D_1) h(D_2) ... h(D_{n-2})\subset K^{n-2}$.
Therefore, $h(\Theta)$ belongs to a fixed compact $K_n:=K^{n-2}$, independent
from the choice of $\Theta$.
\endproof

\hfill

\theorem\label{_BG_qm_Theorem_}
Let $M$ be a closed manifold of strictly negative sectional
curvature,
and $(P, \nabla)$ a principal $G$-bundle with
connection. Fix $x\in M$, and let 
$ q:\;\pi_1(M) \arrow G$
be the non-commutative Barge-Ghys map
associated with $(P, \nabla)$. Then
$q$ is an Ulam quasimorphism (\ref{_Ulam_qm_Definition_}).

\hfill

\proof
Denote by $\tilde M\stackrel \pi \arrow M$ the universal cover of $M$.
Let $a, b\in \pi_1(M)$,
and  $P_a$, $P_b\in \Diff(\tilde M)$ the 
corresponding deck transformations.
Fix a preimage $\tilde x\in \pi^{-1}(x)$,
and denote by $(\tilde P, \tilde \nabla)$
the pullback of $(P, \nabla)$ to $\tilde M$.
By definition, the product $q(ab) q(b)^{-1} q(a)^{-1}$
is represented by the holonomy
of $(\tilde P, \tilde \nabla)$ along the
geodesic simplex connecting
the points $\tilde x, P_a(\tilde x),$ and 
$P_b(P_a(\tilde x))$ in $\tilde M$. 
By \ref{_geodesic_polygon_holonomy_Theorem_}, this quantity belongs to a compact
subset independent from the choice of $x\in M$ and $a, b\in \pi_1(M)$.
\endproof

\hfill

\remark
The set of equivalence classes of non-commutative
Barge-Ghys quasimorphisms is very big.
As follows from
\ref{_Barge_Ghys_classical_Proposition_} below, the usual (``commutative'')
Barge-Ghys quasimorphism is a special case
of the quasimorphism taking values in
a Lie group as defined in \ref{_NC_BG_Definition_}.
By \cite{_Barge_Ghys:bounded_}, \cite{_Calegari:scl_}, 
the vector space spanned by 
commutative Barge-Ghys quasimorphisms up to equivalence is
infinite-dimensional; in \ref{_non_constr_Theorem_}
we construct non-commutative Barge-Ghys quasimorphisms
which cannot be obtained from the commutative ones.

\hfill

\remark 
Using the Riemann-Hilbert correspondence 
(\cite[\S 2.8]{_Ornea_Verbitsky:Book_})
one can associate a flat principal $G$-bundle over $M$ 
with each group homomorphism $\rho:\; \pi_1(M) \arrow G$.
By construction, the holonomy of this flat connection is equal to 
$\rho$ (\cite[\S 2.50]{_Ornea_Verbitsky:Book_}).
Therefore, the corresponding Barge-Ghys quasimorphism is equal to $\rho$.
In other words, all group homomorphisms to a Lie group can be realized
as  Barge-Ghys quasimorphisms.

\subsection{Commutative and
non-commutative Barge-Ghys quasimorphisms }

\definition\label{_classical_BG_Definition_}
Let $M$ be a closed manifold of strictly negative sectional
curvature, and $\theta\in \Lambda^1 M$ a 1-form on $M$.
The {\bf (commutative) Barge-Ghys quasimorphism} 
(\cite{_Barge_Ghys:bounded_})
associated with $\theta$
takes $\gamma \in \pi_1(M)$ to the integral of $\theta$ over
the geodesic path starting and ending at $x$
and homotopic to $\gamma$. 

\hfill

\remark
Note that in all literature
on quasimorphisms, one uses the additive notation:
$|q(xy)- q(x)-q(y)| < C$. For the non-commutative
quasimorphisms, we are forced to use the 
multiplicative notation, \[ q(xy)q(y)^{-1}q(x)^{-1}\in K.\]
When we speak of commutative Barge-Ghys quasimorphisms
in this wider context, this might create a confusion.

\hfill

When $B$ is an oriented real rank-1 vector bundle on $M$,
the commutative Barge-Ghys quasimorphism is actually equal
to the ``non-commutative Barge-Ghys quasimorphism'' defined in
\ref{_NC_BG_Definition_}; we explain this equivalence below.

Topologically, $B$ is always trivial. To trivialize $B$,
we choose a metric on $B$ using the partition of unity,
and trivialize $B$ by taking a positive length-1 section $u$.
Fix the connection $\nabla_0$ on $B$ in such a way that $\nabla_0 u=0$.
Then any connection $\nabla$ on $B$ can be written as
$\nabla= \nabla_0 + \theta$, where $\theta \in \Lambda^1(M)$ is a 1-form.
The holonomy of this connection along a loop $\gamma$ is given
by $\Hol_\gamma(\nabla)= e^{-\int_\gamma \theta}$. 
Indeed, for any section $u_f:=f(t)u$ of $B$ restricted to 
a geodesic segment $\gamma$ parametrized by $t\in [a,b]$,
the equation 
\[ \nabla(u_f)= \frac {df}{dt} u + \theta f u=0
\]
is equivalent to $f' =- f \theta$,
equivalently, $\frac{d\log f}{dt}= - \theta$.

\hfill

\proposition\label{_Barge_Ghys_classical_Proposition_}
Let $(B, \nabla_0)$ be a trivial real rank 1 vector bundle on 
a closed manifold $M$ of strictly negative sectional curvature and
$\theta$ a 1-form on $M$. Consider the connection $\nabla:= \nabla_0 - \theta$.
Then the commutative Barge-Ghys quasimorphism associated
with $\theta$ can be obtained as the logarithm of
the ``non-commutative'' Barge-Ghys quasimorphism associated
with $(B, \nabla)$. 

\proof Indeed, for any loop in $M$, one has
$\Hol_\gamma(\nabla)= e^{-\int_\gamma \theta}$.
\endproof

\hfill

\remark
For any trivialized real rank 1 vector bundle $B$ on $M$,
the connections on $B$ are in bijective correspondence
with 1-forms on $M$. Therefore, 
\ref{_Barge_Ghys_classical_Proposition_}
defines a bijective correspondence between the 
set of (non-commutative) Barge-Ghys quasimorphisms
associated with $B$ and the set of commutative
Barge-Ghys quasimorphisms.

\subsection{Translation length in Gromov hyperbolic groups}

We proceed with a few observations about Gromov hyperbolic groups,
used in the sequel.

\hfill

Let $\Gamma$ be a group
generated by a finite set ${\cal A}$, and 
$C_{\cal A}(\Gamma)$ its Cayley graph. 
{\bf Algebraic translation length} $\tau_{{\cal A}}(\gamma)$ of an
element $\gamma\in \Gamma$ is defined (\cite{_Bridson_Haefliger:book_})
as 
\[
\tau_{{\cal A}}(\gamma):= \lim_{n\to\infty} \frac 1 n d(1_\Gamma, \gamma^n).
\]
The limit exists because the function
$n \mapsto d(1_\Gamma, \gamma^n)$ is subadditive.
It is not hard to see that the map $\gamma\arrow
\tau_{{\cal A}}(\gamma)$ is conjugation invariant 
(\cite[Remark $\Gamma$.3.14 (1)]{_Bridson_Haefliger:book_}).
Further on, we shall use the following result,
also found in \cite{_Bridson_Haefliger:book_}.

\hfill

\proposition\label{_finite_conj_translation_Proposition_}
Let $\Gamma$ be a finitely generated Gromov hyperbolic
group,  and $S_R$
the set of all conjugacy classes of all
$\gamma\in \Gamma$ satisfying $\tau_{{\cal A}}(\gamma) < R$, where $R>0$ is
a real number.
Then $S_R$ is finite.

\proof \cite[Proposition $\Gamma$.3.15]{_Bridson_Haefliger:book_}.
\endproof

\hfill

An element $\gamma\in \Gamma$ is called {\bf primitive}
if it cannot be represented as a power $\gamma= \phi^n$,
for any $n>1$. 

\hfill

\corollary\label{_primitive_Corollary_}
Let $\Gamma$ be a finitely generated Gromov hyperbolic group. 
Then any non-torsion 
element of $\Gamma$ is an integer power
of a primitive element.

\hfill

\proof By \cite[Corollary 3.3.5]{_Calegari:scl_},
$\tau_{{\cal A}}(\gamma)>0$ for all non-torsion $\gamma$.
Clearly, $\tau_{{\cal A}}(\gamma^k)=k\tau_{{\cal A}}(\gamma)$.
By \ref{_finite_conj_translation_Proposition_},
there exists a number $C\in \R^{>0}$ such that
for all non-torsion $u\in \Gamma$ we have $\tau_{{\cal A}}(u) > C$.
Then $\gamma$ cannot be represented as $n$-th power
for any $n > C^{-1}\tau_{{\cal A}}(\gamma)$. Let $m$
be a maximal number such that $\gamma$
is an $m$-th power of an element $\gamma_1\in \Gamma$.
Then $\gamma_1$ is primitive.\footnote{We are grateful to Yves Cornulier
and Sam Nead (\url{https://mathoverflow.net/users/1650/sam-nead})
who suggested a proof of this statement
on Mathoverflow, \cite{_MO:primitive_question_}.}
\endproof

\subsection{HBG-quasimorphism associated with a connection}
\label{_homo_BG_Subsection_}

\definition\label{_homo_Definition_}
A quasimorphism  $q:\; \Gamma \arrow G$ 
is called {\bf homogeneous} if its restriction to any cyclic 
subgroup of $\Gamma$ is a group homomorphism.

\hfill

Let $M$ be a closed Riemannian manifold of strictly negative
sectional curvature, $\Gamma:= \pi_1(M)$. By 
\cite[Corollary 6.2.4]{_Petersen:Riemannian_}, $\Gamma$ is torsion-free.
By \ref{_primitive_Corollary_}, every 
non-torsion element of $\pi_1(M)$ is a power of a
primitive element, which is unique by \cite[Lemma 2.2]{_Bartholdi_Bogopolski_}.


Given a primitive $\gamma \in \pi_1(M)$, let $F_\gamma$
be the shortest free geodesic loop representing
$\gamma$. By \ref{_free_geodesic_unique_Proposition_},
$F_\gamma$ is unique in its free
homotopy class. Clearly, the conjugate elements
of $\pi_1(M)$ correspond to the same free homotopy class.
For each conjugacy class of 
$\gamma$ we fix a choice of
a point $x\in F_\gamma$. When $\gamma= \gamma_1^d$ and
$\gamma_1$ is primitive, we denote by $F_\gamma$
the loop $F_{\gamma_1}$ iterated $d$ times.

Let $(B, \nabla)$ be a bundle with connection. 
Fix a point $p\in M$, and $\gamma\in \pi_1(M)$.
We are going to define a homogeneous quasimorphism
$q:\; \Gamma \arrow GL(B_p)$, where $B_p$ denotes
the fiber of $B$ in $p$, as follows.

Consider a point $x\in F_\gamma$, and let
$\tilde F_\gamma:= \nu_{x, \gamma}\circ F_\gamma\circ \nu^{-1}_{x, \gamma}$
be the 3-segment piecewise geodesic path obtained by connecting $p$ to 
$x$, going around the loop $F_\gamma$ starting and ending
in $x$, and going back to $p$ along $\nu_{x, \gamma}$ in the opposite direction.
Clearly, this path represents $\gamma$ in $\pi_1(M,p)$.
Denote by $q(\gamma)\in GL(B_p)$ the holonomy along
$\tilde F_\gamma$. By construction, $q$ restricted to 
a cyclic subgroup is always a homomorphism.

Note that $q$ depends on the choice of
$x\in F_\gamma$, which has to be fixed for 
each conjugacy class of $\gamma\in \Gamma$.

\hfill

\theorem\label{_homo_Barge_Ghys_Theorem_}
Let $M$ be a closed manifold with strictly negative
sectional curvature, $p\in M$ a base point,
$\Gamma:=\pi_1(M)$, and $q:\; \Gamma
\arrow  GL(B_p)$ the map defined above.
Then $q:\; \Gamma \arrow GL(B_p) $ is an Ulam
quasimorphism. Moreover, $q$ is homogeneous.

\hfill

\begin{figure}[htb]
\centerline{\includegraphics[height=2.8in]{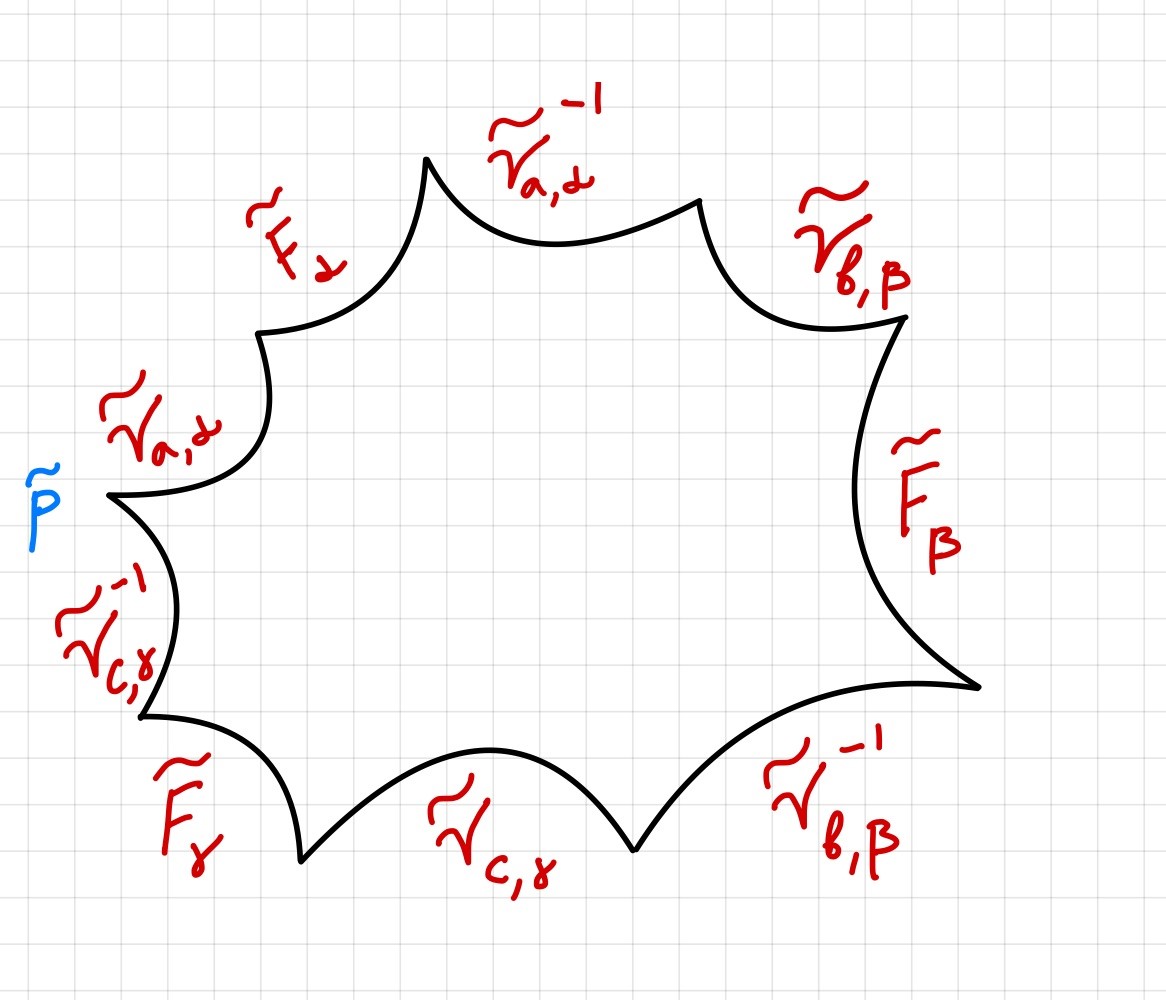}}
\caption{\label{fig:9-gon} Lift of the geodesic polygon to the universal cover.}
\end{figure}

\proof
Let $\alpha, \beta, \gamma=(\alpha\beta)^{-1}$ be elements 
of $\Gamma$. Choose any points $a\in F_\alpha$, $b\in F_\beta$, 
$c\in F_\gamma$. Then $q(\alpha)q(\beta) q(\alpha\beta)^{-1}$
is a holonomy along a contractible geodesic polygon with 9 edges
obtained by going along \[ \nu_{a, \alpha}, F_\alpha, \nu_{a, \alpha}^{-1},
\nu_{b, \beta}, F_\beta, \nu_{b, \beta}^{-1},\nu_{c, \gamma}, F_\gamma, \nu_{c, \gamma}^{-1},
\] 
see Figure \ref{fig:9-gon}.
However, the holonomy along any contractible geodesic
polygon is bounded by
\ref{_geodesic_polygon_holonomy_Theorem_}. 
Homogeneity of $q$ is clear, because 
$q(\gamma^n)$ is the holonomy of $\nabla$ along
the loop $\nu_{x,\gamma}\circ F_\gamma^n\circ \nu_{x,\gamma}^{-1}$.
\endproof

\hfill

\definition\label{_HBG-Definition_}
Let $(B,\nabla)$ be a vector bundle with connection, and
$q$ the Ulam quasimorphism
defined in \ref{_homo_Barge_Ghys_Theorem_}. Then 
$q:\; \Gamma \arrow GL(B_p)$
is called {\bf the HBG-quasimorphism
associated with $(B, \nabla)$.}

\hfill

\remark\label{_G-bundle_HBG_Remark_} 
Instead of vector bundles, \ref{_homo_Barge_Ghys_Theorem_}
can be stated in the setting of principal bundles.
In Section \ref{_principal_Section_} we explained how it is done.
This allows one to modify \ref{_HBG-Definition_}
obtaining a HBG-quasimorphism
associated with a principal $G$-bundle with connection,
for arbitrary Lie group $G$.

\hfill

\remark
A homogeneous quasimorphism
$q:\; \Gamma \arrow \R$
is unique in its equivalence class 
(\cite{_Calegari:scl_}).
Clearly, this is false when $G$
is compact or contains a compact subgroup $K\subset G$:
all homomorphisms from $\Z$ to $K$
are equivalent and homogeneous.
However, the uniqueness is
true  when $G$ is
a simply connected nilpotent Lie group
(\ref{_homo_qm_equivalent=>=_Claim_}).
We expect this to be true for any
simply connected Lie group which has
no non-trivial compact subgroups.

\hfill

\remark
For any homogeneous quasimorphism $q:\; \Gamma \arrow G$
and $r\in G$, the map $r q r^{-1}$ is also  a homogeneous
quasimorphism. The 
HBG-quasi\-mor\-phism 
$q:\; \pi_1(M, p) \arrow G$ associated
with a $G$-bundle depends 
on the choice of the points
$x_i\in F_{\gamma_i}$ on the shortest loop
representing the primitive element $\gamma_i$.
However, this dependence is easy to describe:
for any pair $x_i, y_i \in F_{\gamma_i}$,
consider the geodesic simplex with vertices
$x_i, y_i, p$, with the geodesic from $x_i$ to $y_i$
going in the same direction as in the loop $F_{\gamma_i}$.
Let $r_i$ denote the loop along this simplex, starting
in $p$, going to $x$, then to $y$ and to $p$ again.
Then $q(\gamma_i)$ is the holonomy along
$\nu_{x_i, \gamma_i} F_{\gamma_i} \nu_{x_i, \gamma_i}^{-1}$,
or, equivalently, the holonomy along
$r_i \nu_{y_i, \gamma_i} F_{\gamma_i} \nu_{y_i, \gamma_i}^{-1}r_i^{-1}$.
In other words, replacing the choice of $x_i$ by $y_i$
is equivalent to the conjugation of $q(\gamma_i)$ with the holonomy
along a geodesic simplex $r_i$, which is bounded.
These two quasimorphisms are Ulam equivalent
because they are equivalent to the original Barge-Ghys quasimorphism
(\ref{_HBG_equiv_to_BG_Claim_}). When $G$ is a simply connected
nilpotent Lie group, these HBG quasimorphisms are actually equal
(\ref{_homo_qm_equivalent=>=_Claim_}).

\hfill

In the commutative case, any homogeneous quasimorphism
is conjugation invariant (\cite[\S 2.2.3]{_Calegari:scl_}).
The following result is a non-commutative Barge-Ghys analogue of this
statement.

\hfill

\claim\label{_homo_conjugacy_Claim_}
Let $q:\; \Gamma \arrow G$ be a HBG-quasimorphism, and $\alpha, \beta\in \Gamma$
conjugate elements. Then $q(\alpha)$ is conjugate to
$q(\beta)$.\footnote{Clearly, the conjugating element depends on $\alpha$ and $\beta$.}

\hfill

\proof
Let $F_\gamma$ be the free geodesic loop homotopy
equivalent to $\alpha$ and $\beta$, and 
$x\in F_\gamma$ be the chosen point. Then 
$q(\alpha)$, which is equal to 
holonomy of $\nabla$ along the loop
$\nu_{x, \alpha}\circ F_\gamma\circ \nu_{x,\alpha}^{-1}$,
is conjugate to 
$q(\beta)$, obtained the same way from the loop 
$\nu_{x, \beta}\circ F_\gamma\circ \nu_{x,\beta}^{-1}$.
It is easy to write this conjugation explicitly:
$q(\alpha) = R q(\beta) R^{-1}$, where $R$
is holonomy of $\nabla$ along the loop
$\nu_{x, \alpha} \circ\nu_{x, \beta}^{-1}$.
\endproof

\hfill

When $G\sim \R^n$,
$q(\alpha)= q(\beta)$ for $\alpha, \beta$ conjugate is true by \cite[\S 2.2.3]{_Calegari:scl_}.
However, when $G=S^1$ (or other compact group),
it is easy to construct a homogeneous Ulam quasimorphism
$q:\; \Gamma \arrow G$ which does not have this property. For example,
let $\Gamma= {\Bbb F}_2$ be a free group, and define
an arbitrary map $q:\; {\goth P} \arrow S^1$, where ${\goth P}$ is the
set of all primitive words in ${\Bbb F}_2$.
If $q(W)= q(W^{-1})^{-1}$, this map extends to 
a homogeneous Ulam quasimorphism, and all
homogeneous Ulam quasimorphisms $q:\; {\Bbb F}_2 \arrow S^1$, are obtained this way.
Clearly, this map  is not necessarily invariant under the conjugation.
This implies, in particular, that \ref{_homo_conjugacy_Claim_} is
false for a general Ulam quasimorphism, if $G$ has a
compact subgroup.

\hfill 

\question
Can \ref{_homo_conjugacy_Claim_} be generalized to 
an arbitrary homogeneous Ulam quasimorphism
 $q:\; \Gamma \arrow G$, when $G$ is a Lie group
which has no compact subgroups?

\hfill

Any commutative quasimorphism is equivalent to 
its homogenization, which is unique in its equivalence
class (\cite[\S 2.2.3]{_Calegari:scl_}). For Ulam quasimorphisms taking values in
a non-commutative group, the homogenization procedure
is not yet known. However, the Barge-Ghys quasimorphisms
admit a ``homogenization'', which is unique in its equivalence 
class by \ref{_homo_qm_equivalent=>=_Claim_} below.

\hfill

\claim\label{_HBG_equiv_to_BG_Claim_}
Let $(B,\nabla)$ be a bundle (vector or a principal $G$-bundle) 
with connection over a closed manifold of strictly negative sectional curvature,
 $\Gamma=\pi_1(M,p)$, and  $q_\nabla:\; \Gamma\arrow G$ the corresponding  Barge-Ghys quasimorphism.
After making the relevant choices, we obtain the 
HBG-quasimorphism $q:\; \Gamma\arrow G$. Then $q$ is equivalent
to $q_\nabla$, and, moreover, 
$q(\gamma) q_\nabla(\gamma)^{-1}\in K_4$, where $K_4\subset G$ is the compact
subset defined in \ref{_geodesic_polygon_holonomy_Theorem_},
which depends only on $M$, $B$ and $\nabla$.

\hfill

\proof
Let $\gamma\in \Gamma$ be any element,
which is represented as $\gamma=\gamma_0^n$, where
$\gamma_0$ is primitive, and let $F_{\gamma_0}$
be the free geodesic representing the same free homotopy class.
Fix a geodesic segment $\nu_{\gamma_0,x}$ connecting $p\in M$ to
a point $x$ in $F_{\gamma_0}$, and
let $\tilde \gamma$ the geodesic segment
connecting $p$ to itself and homotopic to $\gamma$.
Then $q_\nabla(\gamma)$ is a holonomy of $\nabla$ along
$\tilde \gamma$, and $q(\gamma)$ is the holonomy
of $\nabla$ along the geodesic chain 
$\nu_{\gamma_0,x}\circ F_{\gamma_0}^n\circ \nu_{\gamma_0,x}^{-1}$.
Geometrically, the segment $F_{\gamma_0}^n$ 
is one geodesic interval, and therefore the
polygon $\nu_{\gamma_0,x}\circ F_{\gamma_0}^n\circ \nu_{\gamma_0,x}^{-1} \circ \tilde\gamma^{-1}$
has 4 geodesic sides. By \ref{_geodesic_polygon_holonomy_Theorem_},
the holonomy along the sides this polygon belongs to $K_4$. 
\endproof

\hfill

\remark 
For homogeneous quasimorphisms taking value in $\R$ or $\R^n$
(elsewhere, we call them ``commutative''),
equivalence implies equality: if $q_1, q_2:\; \Gamma \arrow \R^n$
are equivalent homogeneous quasimorphisms, they satisfy $q_1=q_2$,
because
\[ q_1(x)- q_2(x) = \lim_{n\to \infty}
\frac{q_1(x^n)-q_2(x^n)}{n}=0.
\] 
The same is true for all homogeneous Ulam quasimorphisms
$q:\; \Gamma \arrow G$, if $G$ is nilpotent and simply connected:

\hfill

\claim\label{_homo_qm_equivalent=>=_Claim_}
Let $q_1, q_2:\; \Gamma \arrow G$
be homogeneous quasimorphisms which are  
equivalent. Assume that $G$ is a simply connected nilpotent Lie group.
Then $q_1=q_2$.\footnote{We expect that this statement
is true whenever $G$ is a Lie group such that the
exponential map $\exp:\; \g \arrow G$  is a homeomorphism.}

\hfill

\pstep
Let $K\subset G$ be a compact subset such that  
$q_1(x) q_2(x)^{-1} \in K$ for all $x\in \Gamma$.
Since $q_i(x^n)= q_i(x)^n$, this gives
 $q_1(x)^n q_2(x)^{-n}\in K$. Then 
\ref{_homo_qm_equivalent=>=_Claim_}
would follow if we prove that for 
any distinct $a, b \in G$, the sequence
$a^n b^{-n}$ is unbounded, that is, does not
belong to a compact set.

%

\hfill

{\bf Step 2:}
Let $G_i$ be the lower central series for $G$.
Applying induction and using 
the long exact sequence 
\[  \arrow \pi_2(G_i/G_{i-1}) \arrow \pi_1(G_{i-1}) \arrow \pi_1(G_i) \arrow \pi_1(G_i/G_{i-1}) \arrow 0
\]
for the Serre's 
fibration $G_{i-1} \arrow G_i \arrow G_i/G_{i-1},$
we immediately obtain that the groups
$G_i$ and $G_i/G_{i-1}$ are simply connected for all $i$.
Since $G_i/G_{i-1}$ is abelian, this group is isomorphic to $\R^n$.

\hfill

{\bf Step 3:}
We  use the induction on $\dim G$.
When $G$ is 1-dimensional,
this group is commutative, giving $G=\R$;
if the sequence $a^n b^{-n}\in \R$ is bounded,
we have $a=b$.

Let $Z\subset G$ be the center; it is
non-trivial, because $G$ is nilpotent. Since
$a^n b^{-n}$ is bounded modulo $Z$, we may apply
the induction hypothesis to the representatives
of $a, b$ in $G/Z$ and obtain that $a=b$ modulo $Z$.
Then $az=b$, where $z\in Z$, hence the sequence 
$a^n b^{-n}= z^n$ is bounded. Since $Z=\R^n$ (Step 2), this
may happen only if $z$ is trivial. Then $a=b$.
\endproof

\hfill

\remark
Let $E$ be a principal $G$-bundle with connection with $G=\R^{>0}$, and
$B$ the associated line bundle with induced connection. Then the usual
(commutative) Barge-Ghys quasimorphism $q$, taking values in
the additive group $\R$, is obtained from our Barge-Ghys quasimorphism
$q_\nabla$ by taking the logarithm. 
Let $q_0$ be the the corresponding HBG-quasimorphism.
Then the logarithm $\log q_0$
is homogeneous and equivalent to $q$ by \ref{_HBG_equiv_to_BG_Claim_}.
Since a homogeneous commutative quasimorphism is unique in its 
equivalence class (\cite[Lemma 2.2.1]{_Calegari:scl_}),
this implies that $\log q_{0}$ is equal to the homogenization of $q$.

\section{Constructible Ulam quasimorphisms}

In the sequel, when we speak of nilpotent Lie groups,
we always mean connected, simply connected, algebraic nilpotent Lie
groups over $\R$. Sometimes we emphasise ``algebraic'', but
this is not that necessary: every connected simply connected nilpotent Lie group
admits a unique algebraic structure by 
\cite[Theorem on page 12]{_Hochschild:algebraic_Lie_}.

\subsection{Zariski closure of discrete subgroups}

\definition
Let $Y$ be a real algebraic variety, and 
$X \subset Y$ a subset. {\bf The Zariski closure}
of $X$ in $Y$ is the intersection of all
real algebraic subvarieties $X_i\subset Y$
containing $X$. A subset $X \subset Y$
is {\bf Zariski dense} if its Zariski 
closure is $Y$.

\hfill

\remark
Let $\Gamma \subset G$ be a subgroup of a real algebraic
group. Then its Zariski closure $\bar \Gamma$ is an algebraic subgroup
of $G$. Indeed, the group laws put algebraic constraints
on $\bar \Gamma$, hence the smallest algebraic subvariety
of $G$ containing $\Gamma$ is closed under the  group
operations.

\hfill

The following notion is going to be used in 
the proof of \ref{_non_constr_Theorem_} below.
For the readers' convenience, we recall the definition
of Hausdorff distance (\cite{_Gromov:Riemannian_}).

\hfill

\definition\label{_Hausdorff_distance_Definition_}
Let $Z_1, Z_2$ be subsets of a metric space $M$.
Denote by $Z_i(\epsilon)$ the $\epsilon$-neighbourhood of $Z_i$, 
and let 
\[ d_H(Z_1, Z_2):= \inf\{ \epsilon\in \R^{>0}\cup \infty
\ \ |\ \ Z_1(\epsilon)\supset Z_2 \text{\ and \ }
Z_2(\epsilon)\supset Z_1
\}.
\]
The number $d_H(Z_1, Z_2)$ is called {\bf the Hausdorff distance}
between $Z_1$ and $Z_2$. It is not hard to see that
$d_H$ defines a metric, taking values in $\R^{>0}\cup \infty$,
on the set of all closed subsets of $M$.

\hfill

\definition
Two subsets $X, Y$ of a metric
space $M$ are {\bf coarse equivalent} if 
the Hausdorff distance $d_H(X,Y)$ is finite.
When $M$ is a Lie group, and $d$ a
left-invariant Riemannian distance, this is equivalent
to $K\cdot X \supset Y$ and $K \cdot Y \supset X$
for a compact subset $K \subset G$. 

\hfill

\example
Any two rank 2 discrete lattices in $\R^2$ with the usual 
Euclidean metric are coarse
equivalent. However, $\R^2$ is not coarse equivalent to a point.

\hfill

\example More generally, $\R^k$ taken with the standard
Euclidean metric is not coarse equivalent to
any subspace $\R^n \subsetneq \R^k$.

\hfill

We need to relax
the notion of coarse equivalence
slightly to accomodate the metrics which are
not bi-invariant.

\hfill

\definition
Let $G$ be a Lie group. We say that
$X \subset G$ is {\bf bi-coarse equivalent}
to $Y\subset G$ if there exists a compact subset
$K \subset G$ such that $K\cdot X \cdot K \supset Y$
and $K\cdot Y \cdot K \supset X$.

\hfill

\definition\label{_finite_dist_Zariski_Definition_}
Let $G$ be a real algebraic Lie group,
and $\Lambda\subset G$ a Zariski dense subgroup.
We say that $\Lambda$ is {\bf bi-coarse Zariski dense in $G$}
if the following property holds.
If a subset $\Lambda' \subset G$ is
bi-coarse equivalent to $\Lambda$, then $\Lambda'$ is also 
Zariski dense in $G$.

\hfill

We prove that a lattice in a nilpotent Lie group is bi-coarse
Zariski dense. Note that simply connected nilpotent Lie
groups are algebraic, with the algebraic structure induced
by the standard algebraic structure on its Lie algebra.
Indeed, for a simply connected nilpotent Lie group, the exponential map
from its Lie algebra to the Lie group
is polynomial and invertible, and the 
inverse map, called the logarithm, 
is also polynomial (\cite[Proposition 1.2.8]{_Corwin_Greenleaf_}).

\hfill

\definition
A subset $\Lambda \subset G$  of a Lie group is called {\bf
  bi-cocompact}
if there exists a compact subset $K \subset G$ such that
$K \cdot \Lambda\cdot K = G$.

\hfill

\remark 
Further on, we use the following elementary observation.
Clearly, the subgroup $[G,G]\subset G$ is normal, and
for any $H\subset G$ the group $H\cdot[G,G]$ generated
by $H$ and $[G,G]$ is also normal. Indeed, any subgroup
$G_1\subset G$ which contains $[G, G]$ is normal, because
$x^{-1} yxy^{-1}\in [G,G]\subset G_1$ for any $x\in G_1$ and $y\in G$,
and therefore $yxy^{-1}$ also belongs to $G_1$.

\hfill

\claim\label{_non_cocompact_subgroup_nilpo_Claim_}
Let $G$ be a connected, simply connected algebraic nilpotent
Lie group, and $H\subset G$ a proper algebraic subgroup.
Consider the minimal Lie subgroup $H\cdot[G,G]$ containing $H$ and $[G,G]$.
Then $G/(H\cdot[G,G])\cong \R^n$ for some $n>0$. 

\hfill

\proof 
A quotient of a connected, simply connected 
nilpotent Lie group by an algebraic subgroup
is connected, simply connected, because this
algebraic subgroup is connected and simply connected
(\cite[Proposition 1.2.8]{_Corwin_Greenleaf_}).
Since $G/(H\cdot[G,G])$ is commutative,
it is isomorphic to $\R^n$. It remains
to show that $n>0$. Otherwise
the projection of $\goth h$ to $\frac{\g}{[\g, \g]}$ is 
surjective, where $\g = \Lie(G)$ and $\goth h = \Lie(H)$.
However, any set of elements generating $\frac{\g}{[\g, \g]}$ 
also generates $\g$ (\cite{_Koch:generators_}), hence in this case $G =H$,
which is a contradiction because $H\subset G$ is a proper subgroup.
\endproof

\hfill

\remark
Let $H\subset G$ be a Zariski closed
subgroup in a connected, simply connected nilpotent
algebraic Lie group. Then $H$ is also connected
and simply connected (\cite[Proposition 1.2.2]{_Corwin_Greenleaf_}).

\hfill

We need the following trivial sublemma.

\hfill

\sublemma \label{_neigh_of_supspace_Sublemma_}
Consider vector spaces $\R^k\subsetneq \R^n$,
and let $K\subset \R^n$ be a compact subset.
Then $K + \R^k \neq \R^n$.
\endproof

\hfill

\lemma\label{_cocompact_hence_zar_dense_Lemma_}
Let $\Lambda\subset G$ be a bi-cocompact subset in a simply
connected  nilpotent Lie group $G$. Then the group
generated by $\Lambda$ is Zariski dense. 

\hfill

\proof
Let $H\subset G$ be the Zariski closure of the group
generated by $\Lambda$. We need to show that $H= G$.
Since $H$ is bi-cocompact, $K \cdot H \cdot K=G$ for some compact
subset $K \subset G$. Consider the projection map 
$\pi:\; G \arrow G/[G, G]$. Let $K_1:= \pi(K)$; this set is
compact because $K$ is compact and $\pi$ is continuous.
The group $\pi(H)$, which is a commutative
nilpotent Lie group, is isomorphic to $\R^n$. Assume, on contrary, that
$H\neq G$.
By \ref{_non_cocompact_subgroup_nilpo_Claim_}, the image of
$H$ in $G/[G, G]$ is strictly smaller than $\pi(G)$, which is also 
homeomorphic to $\R^m$, for some $m>n$.
Now, $K_1 \pi(H) K_1 \neq \pi(G)$ follows from
\ref{_neigh_of_supspace_Sublemma_}. This gives a contradiction
with $K \cdot H \cdot K=G$.
\endproof

\hfill

\corollary\label{_bi-coarse_dense_nilp_Corollary_}
Let $G$ be a simply connected nilpotent Lie group,
and $\Gamma \subset G$ a lattice, which is 
 cocompact (\cite[Theorem 2.2.6]{_Wang:597_}). Then 
$\Gamma$ is bi-coarse Zariski dense.

\hfill

\proof
Since $\Gamma$ is cocompact, it is bi-cocompact.
Any subset  $S\subset G$ which is bi-coarse equivalent
to a bi-cocompact set is bi-cocompact. Therefore,
$S$ is Zariski dense by 
\ref{_cocompact_hence_zar_dense_Lemma_}.
This implies that $\Gamma$ is bi-coarse Zariski dense.
\endproof

\hfill

\subsection{Virtually conjugation equivalent cyclic subgroups}

In the sequel, we need the following notion.

\hfill

\definition
Let $A, B\subset G$ be infinite cyclic subgroups
of a group $G$. We say that $A$ is
{\bf virtually conjugation equivalent (VCE equivalent) to $B$}, denoted
$A \sim_{\vce} B$,
if there exist $u\in G$ such that the intersection
$A \cap B^u$ is infinite, where $B^u$ denotes
the subgroup obtained from $B$ by conjugation with $u$.
We also write $x\sim_\vce y$ for elements $x, y \in G$
when $x, y$ generate infinite cyclic subgroups which
are VCE equivalent.

\hfill

\remark\label{_commens_elements_Remark_}
In \cite{_DGO:hyperb_embedded_Memoirs_},
two elements of a group are called {\bf commensurable}
if non-zero powers of these elements are conjugate.
This is equivalent to VCE equivalence of cyclic subgroups generated
by these elements.

\hfill

\claim
The relation $\sim_{\vce}$ is, indeed, an equivalence relation.

\proof
Let $A, B, C$ be three cyclic subgroups which satisfy 
$A \sim_{\vce} B$ and $B \sim_{\vce} C$. 
Then $A$ is commensurable
with $B^u$ and  $B$ is commensurable
with $C^v$, that is, the intersections $A \cap B^u$ 
and $B\cap C^v$ are infinite. 
Note that any infinite subgroup in a cyclic group
has finite index, hence to show that $A\sim_{\vce}C$
it would suffice to show that $A \cap C^{vu}$ is infinite.

By definion, $C^{vu}\cap B^u$ 
is a finite index subgroup in $B^u$. On the other hand,
$B^u$ contains $A \cap B^u$
as a finite index subgroup. 
Since an intersection of two infinite subgroups
in $\Z$ is always infinite, $A\cap B^u\cap C^{vu}$
is infinite, proving the claim.
\endproof

\hfill

\remark
If $\rk H^1(\Gamma, \Q)\geq 2$,
it is easy to find infinitely many 
cyclic subgroups which are pairwise VCE 
non-equivalent. Indeed, if the line
in $H^1(\Gamma, \Q)$ generated
by a cyclic subgroup $A$ is not 
collinear with the line generated by $B$,
this implies that $A \not\sim_\vce B$.
To extend this statement to more general
hyperbolic groups, we need the following
argument.

\hfill

The notion of hyperbolic embedding was defined in \cite{_Dahmani_Guirardel_Osin_}.
The following result is one of the applications of 
hyperbolic embeddings

\hfill

\claim \label{_extending_quasimorphisms_Claim_}
Let $H$ be hyperbolically embedded subgroup of $\Gamma$,
then every (commutative) quasimorphism on $H$ extends to a quasimorphism on $\Gamma$.

\proof
The  quasimorphisms from $H$ to $\R$ are controlled by kernel 
of the map $H^2_b(H, \R)\arrow H^2(H, \R)$, called {\bf the exact
reduced bounded cohomology}, which is clear from \eqref{_bounded_exact_Equation_}.
Every class in second exact reduced bounded cohomology of $H$ can be extended
to $H^2_b(\Gamma, \R)$, as follows from \cite{_Hull_Osin_}.
This result was generalized to all cohomology groups in
\cite{_Frigerio_Pozzetti_Sisto_}. \endproof

\hfill

In \cite{_Osin_}, D. Osin defined acylindrically hyperbolic groups, 
which includes all hyperbolic groups which are not virtually cyclic.
The fundamental group of a closed manifold of strictly negative
sectional curvature, which is by convention assumed of dimension
$>1$, is hyperbolic and
not virtually cyclic, hence it is acylindrically hyperbolic.
In \cite[Theorem 1.2]{_Osin_} and \cite[Theorem 2.24]{_Dahmani_Guirardel_Osin_},
it was shown that any acylindrically hyperbolic group $\Gamma$
contains a hyperbolically embedded subgroup 
$H={\Bbb F}_2\times K$, where $K$ is some finite group.
This brings the following theorem.

\hfill

\theorem\label{_free_group_extending_quasimorphisms_Theorem_}
Let $\Gamma$ be an  acylindrically hyperbolic group 
(this includes fundamental groups of a closed manifold
of strictly negative sectional curvature.)
Then there exists
a free subgroup $H={\Bbb F}_2\subset \Gamma$
such that any quasimorphism on $H$ can be
extended to a quasimorphism on $\Gamma$.
\endproof

\hfill

\remark
Brooks quasimorphisms are quasimorphisms
on a free group, defined in \cite{_Brooks:bounded_};
see \cite{_Calegari:scl_} for more details.
We use \ref{_free_group_extending_quasimorphisms_Theorem_}
only Brooks quasimorphisms; however, for Brooks
quasimorphisms this result is implicit
already in  Bestvina-Fujiwara
\cite{_Bestvina_Fujiwara_}, 
where the free group is a Schottky subgroup.

\hfill

The following result immediately follows
from the definition of Brooks quasimorphism.

\hfill

\claim\label{_Brooks_not_subwords_Claim_}
Let $W_1$ and $W_2$ be two reduced words, considered
as elements of ${\Bbb F}_2$.
Denote by $q_1:\; {\Bbb F}_2\arrow \R$ the homogenization of the 
Brooks quasimorphism associated
with $W_1$.  Assume that $W_1$ is
not a subword of $W_2^i$, for all $i\in \Z$.
Then $q_1(W_2) =0$.

\proof The proof is left as an exercise to the reader.
\endproof

\hfill

This implies the following proposition which will be used in the sequel.

\hfill

\proposition \label{_a^nb^m_non_conj_Proposition_}
Let $\Gamma$ be an  acylindrically hyperbolic group.\footnote{
This includes fundamental groups of a closed manifold
of strictly negative sectional curvature.} 
Let $k, k_1, l, l_1$ be non-negative integers
such that neither $k=0, k_1=0$ nor $l=0, l_1=0$,
nor $k=k_1, l=l_1$ holds.
Then there exists
a free subgroup $H={\Bbb F}_2\subset \Gamma$
generated by $a, b\in \Gamma$ such that
 $a^kb^l \not \sim_\vce a^{k_1} b^{l_1}$,
where the relation $\not\sim_\vce$ is taken 
in $\Gamma$.

\hfill

\proof
By \ref{_free_group_extending_quasimorphisms_Theorem_},
there exists a subgroup ${\Bbb F}_2\subset \Gamma$
such that any quasimorphism is extended from 
${\Bbb F}_2$ to $\Gamma$. Let $W_1=a^{k_1} b^{l_1}$ and
$W_2=a^kb^l$. The negation of the condition ``$k=0, k_1=0$ or $l=0, l_1=0$,
or $k=k_1, l=l_1$'' is equivalent to ``$W_1$
is not a subword of $W_2^i,$ and $W_2$
is not a subword of $W_1^i$, for all $i\in \Z$.''
Let $q_1, q_2:\; {\Bbb F}_2\arrow \R$ be the
homogenizations of the corresponding Brooks quasimorphisms.
By construction, $q_i$ can be extended
to quasimorphisms  $q_1, q_2:\; \Gamma\arrow \R$.
Since homogeneous quasimorphisms are conjugate
invariant, and $q_i(W_j) = \delta_{ij}$,
the cyclic group $\langle W_1\rangle$
intersects trivially  with any conjugate to 
$\langle W_2\rangle$.
\endproof

\hfill

As an immediate corollary, we obtain the following
assertion, which is also implied by 
\cite{_Bestvina_Fujiwara_}.\footnote{For non-elementary
hyperbolic groups, this result also follows from 
\cite{_Epstein_Fujiwara_}.}

\hfill

\corollary\label{_conj_classes_infinitely_many_Corollary_}
Let $\Gamma$ be an acylindrical hyperbolic group.
Then there exists infinitely many cyclic subgroups
$U_i$ which satisfy $U_i \not \sim_\vce U_j$ for any $i\neq j$.
\endproof

\hfill

\remark
In \cite[Corollary 6.12]{_DGO:hyperb_embedded_Memoirs_},
a similar result is shown: the authors construct
an arbitrarily large collection of elements in a hyperbolic group
which are non-commensurable in the sense of 
\ref{_commens_elements_Remark_}.

\subsection{Non-constructible Ulam quasimorphisms}

In
\cite{_Fujiwara_Kapovich_}, Kapovich and Fujiwara
define a constructible quasimorphism taking
values in a discrete group. We extend their
definition to any Lie group.

\hfill

\definition\label{_constructible_Definition_}
Let $q:\; \Gamma \arrow G$ be an Ulam quasimorphism. 
Consider a finite index subgroup
$\Gamma_0\subset \Gamma$, and let $H \supset q(\Gamma_0)$
be a subgroup of $G$ containing $q(\Gamma_0)$.
Assume that $H$ contains a normal subgroup $A$,
abelian and central in $H$, such that the composition 
$\Gamma_0 \stackrel q \arrow H \arrow H/A$
is equivalent to a homomorphism $\Gamma_0 \arrow H/A$.
Then $q$ is called {\bf a constructible quasimorphism}.

\hfill

An element $a\in \Gamma$ of a group is called {\bf primitive}
if for any $n\in \Z$ and $b\in \Gamma$ such that $b^n=a$, one has $n=\pm 1$.
A quasimorphism $q:\; \Gamma \arrow G$
is called {\bf homogeneous} if its restriction to any
cyclic subgroup is a homomorphism (\ref{_homo_Definition_}). 
We gave a detailed treatment of homogeneous quasimorphisms in Subsection
\ref{_homo_BG_Subsection_}, where we defined a new type of
quasimorphisms, called HBG-quasimorphism, for ``homogeneous
Barge-Ghys''. We use the acronym HBG to avoid the confusion, because
the HBG-quasimorphisms are not, in fact,
``Barge-Ghys quasimorphisms'', 
in the sense of \ref{_NC_BG_Definition_}. One should think
of a HBG-quasimorphism as of a ``homogenization'' of a Barge-Ghys 
quasimorphism. However, the notion of homogenization, which is
well known for the usual ($\R$-valued) quasimorphisms, is not
defined (yet) for the quasimorphisms taking values in a Lie group.

In Subsection \ref{_constructing_Subsection_},
we prove the following theorem.

\hfill

\theorem\label{_existence_Ulam_BG_Theorem_}
Let $M$ be a closed manifold of strictly negative
curvature, $G$ a connected non-abelian Lie group, and 
$x_1, ..., x_n\in \Gamma:= \pi_1(M)$
a collection of primitive elements satisfying $x_i \not \sim_\vce x_j$
for $i\neq j$.\footnote{Since $\dim M>1$, there
exist infinitely many such subgroups, 
as follows from \ref{_conj_classes_infinitely_many_Corollary_}.}
Take any collection of elements $g_i\in G$, with $i= 1, ..., n$.
Then 
\begin{description}
\item[(i)]
There exists a connection $\nabla$ on a trivial 
principal bundle $P$ such that the corresponding 
HBG-quasimorphism $q_\nabla:\; \Gamma \arrow G$
takes $x_i$ to $g_i$, $i=1, 2, ..., n$.
\item[(ii)]
Moreover, for any countable family of 
elements $z_i$ such that $z_i \not \sim_\vce x_j$
for all $z_i, x_j$, the connection $\nabla$ can be chosen
such that $q_\nabla(z_i)$ is not central in $G$.
\end{description}

\proof See Subsection \ref{_constructing_Subsection_} 
below. \endproof

\hfill

We use this result to prove \ref{_non_constr_Theorem_} below.

\hfill

\definition
A real nilpotent
Lie group is called {\bf rational} if its Lie algebra $\g$ contains
a rational Lie subalgebra $\g_\Q$ such that $\g=\g_\Q\otimes_\Q \R$.
By Maltsev's theorem, this property is equivalent to existence
of cocompact lattices (\cite{_Corwin_Greenleaf_}). 

\hfill

\example
Consider the {\bf Heisenberg group} ${\cal H}$,
that is, the group of upper triangular matrices
$3\times 3$. This is a 3-dimensional nilpotent
Lie group. The subgroup ${\cal H}_\Z$ of integer
upper triangular matrices is cocompact in ${\cal H}$,
which is an exercise left to the reader.
This implies, in particular, that  ${\cal H}$ is rational.

\hfill

\theorem \label{_non_constr_Theorem_}
Let $G$ be a simply connected, connected, non-abelian
rational real nilpotent Lie group, and 
$\Gamma:= \pi_1(M)$,  where $M$ is a closed
manifold of strictly negative sectional curvature, $\dim_\R M >1$.
Then there exists a non-constructible HBG-quasimorphism
$q_\nabla:\; \Gamma \arrow G$.

\hfill

\pstep
By Maltsev's theorem, $G$ contains a lattice $\Lambda$
(\cite[Theorem 5.8.1]{_Corwin_Greenleaf_}),
which is a posteriori cocompact (\cite[Corollary
  5.4.6]{_Corwin_Greenleaf_}). This lattice is 
finitely generated (\cite[Corollary
  5.1.7]{_Corwin_Greenleaf_}). Denote its
generators by $g_i$.  By
\ref{_conj_classes_infinitely_many_Corollary_},
$\Gamma$ contains infinitely many $x_i$ 
which satisfy $x_i \not\sim_\vce x_j$ for all $i\neq j$. Using
\ref{_existence_Ulam_BG_Theorem_}, we find
$x_1, ..., x_n\in \Gamma$ and 
a HBG-quasimorphism 
$q_\nabla:\; \Gamma \arrow G$
taking each $x_i$ to $g_i$.
We choose, in addition, elements $a, b \in \Gamma$
such that $a^l b^m\not \sim_\vce a^{l_1} b^{m_1}$
unless $l=l_1, m=m_1$ or  $l=l_1=0$ or $m=m_1=0$  
(\ref{_a^nb^m_non_conj_Proposition_}). 
Since the choice of $x_i$ is arbitrary,
we may assume also that $x_i \not \sim_\vce a$ and
$x_i \not \sim_\vce b$ for all $i$.

Denote the identity element in $G$ by $1_G$.
We apply
\ref{_existence_Ulam_BG_Theorem_}
to $x_i$ as above, $z_i$ equal to
 $a^l b^m$ such that $l, m \neq 0$. Then 
we obtain a HBG-quasimorphism such that $q_\nabla(x_i)=g_i$,
$q_\nabla(a)=1_G$, $q_\nabla(b)=1_G$,
and $q_\nabla(a^ib^j)$ is not central for 
infinitely many $i, j>0$. 

We are going to show that $q_\nabla$
is not constructible.

\hfill

{\bf Step 2:}
By contradiction, suppose that $q_\nabla$ is constructible, that is,
there exists a finite index subgroup 
$\Gamma_0\subset \Gamma$
and
$q_\nabla\restrict{\Gamma_0}$ 
satisfies the assumptions of \ref{_constructible_Definition_}.

Since $\Gamma_0$ is of finite index 
in $\Gamma$, it contains powers $x_i^{n_i}$
of each $x_i$, for some $n_i \in  \Z^{>0}$. 
Using the central series for the lattice $\Lambda\subset G$,
we can easily see that the elements $q_\nabla(x_i^{n_i})= g_i^{n_i}$
generate a finite index sublattice $\Lambda_0\subset \Lambda$, which 
is bi-coarse equivalent to $G$ because it is cocompact.

Since $\Lambda_0$ is bi-coarse Zariski dense 
(\ref{_bi-coarse_dense_nilp_Corollary_}), 
any map $\Gamma_0 \arrow G$ which is 
equivalent to $q_\nabla$ 
has a Zariski dense image. 

Since $q_\nabla$ is constructible, there exists
a homomorphism $h:\; \Gamma_0 \arrow H/A$, where 
$H\subset G$ is a subgroup and $A$ its central subgroup,
which can be lifted to a homomorphism $h_1:\; \Gamma_0\arrow H$
which is Ulam equivalent to $q_\nabla\restrict {{\Gamma_0}}$.
Then $\im h$ is Zariski
dense in $G$, as indicated above. Therefore
the group $H$ is also Zariski dense in $G$,
and $A$ belongs to the center $Z\subset G$.

Clearly,
the composition of $h:\; \Gamma_0 \arrow G/A$
and  $\pi:\; G/A \arrow G/Z$ is a homomorphism.
It remains to show that the composition of
$q_\nabla\restrict {\Gamma_0}$ and  $\pi:\; G \arrow G/Z$
is not equivalent to any homomorphism,
which will give a contradiction.

\hfill

{\bf Step 3:} A sequence $\{x_i\in G\}$ is called {\bf bounded}
if it belongs to a compact subset of $G$. For any compact subset 
$K \subset G$ in a topological group $G$,
and any sequences $\{x_i\}$, $\{y_i \in K x_i K\}$,
the sequence $\{x_i\}$ is bounded if and only if $\{y_i\}$
is bounded. Indeed, $y_i \in  K x_i K$ $\Leftrightarrow$
$x_i \in  K^{-1} y_i K^{-1}$, and for any bounded $\{x_i\in K_1\}$,
the sequence $y_i \in K x_i K\subset K K_1 K$ is clearly bounded.

\hfill

{\bf Step 4:} A simply connected nilpotent Lie group $G$ does not
contain a bounded subgroup. Indeed, the closure of a bounded
group is compact, because a closed bounded set is compact.
By Cartan's theorem, its closure is a Lie subgroup $K \subset G$.
Let $G \supset G_1 \supset G_2 \supset ...$ be the lower
central series for $G$, with $G_{i+1} = [G_i, G]$, and $i$ the smallest number such that 
$K \subset G_i$. Then the projection of $K$ to $\frac{G_i}{[G_i, G_i]}$
is a non-trivial compact subgroup. However, $\frac{G_i}{[G_i, G_i]}=\R^n$,
and $\R^n$ does not have non-trivial compact subgroups.

\hfill

{\bf Step 5:}
Let $q_1:\; \Gamma \arrow G/Z$ denote 
the composition of $q_\nabla:\; \Gamma \arrow G$ and the projection to $G/Z$.
Arguing by contradiction, suppose that $h:\; \Gamma_0 \arrow G/Z$ is a homomorphism 
which is equivalent to $q_1$, where $\Gamma_0\subset\Gamma$ is a finite index subgroup.
Then there exists a compact $K\subset G$ such that
\begin{equation}\label{_AQM-equiv_Equation_}
h(z) \in K q_1(z) K
\end{equation} 
for all $z\in \Gamma_0$.
Since $\Gamma_0$ has finite index in $\Gamma$,
there exist $i, j>0$ such that $a^i, b^j\in \Gamma_0$.
Since $q_1(a^{ki})=q_1(b^{kj})=1_G$ for
all $k\in \Z$,
this gives $h(a^{ik})=h(a^i)^k \in K\cdot K$ and 
$h(b^{jk})=h(b^j)^k \in K\cdot K$.
However, for any element $x$ in 
a simply connected nilpotent Lie group, a sequence
$\{x^n, n\in \Z\}$
can be bounded only if $x=1_G$, hence
$h(a^i)=1_G$ and $h(b^j)=1_G$. Since $h$ is a homomorphism,
this implies that
$h(a^i b^j)=1_G$. By construction, the HBG-homomorphisms
are ``homogeneous'', that is, satisfy $q_\nabla(x^n)=(q_\nabla(x))^n$.
This implies that $q_1$ is also homogeneous.
Since $q_1 (a^i b^j)\neq 1_G$, 
the sequence $q_1((a^i b^j)^n), n\in \Z$ is not
bounded (Step 4), which implies that the sequence
$h((a^i b^j)^n) \in K q_1((a^i b^j)^n)K$
is also unbounded (Step 3), giving a contradiction.
\endproof

\hfill

\remark \label{_algebra_equiv_Remark_}
In \ref{_non_constr_Theorem_}, we construct 
a HBG-quasimorphism which is
not equivalent to a constructible one. However,
the proof we use brings a more powerful result,
which was the chief aim of an earlier version 
of this paper. Let $q, q':\; \Gamma \arrow G$ be
Ulam quasimorphisms. We say that $q, q'$ are
{\bf algebraically equivalent}
if there exists a compact subset $K\subset G$
such that for all $x\in \Gamma$, one has
$q(x) \in K q'(x) K$. In \ref{_non_constr_Theorem_},
we construct a HBG-quasimorphism
$q_\nabla$ which is not algebraically equivalent to
a constructible quasimorphism.

\subsection{Different notions of quasimorphisms:
geometric, algebraic and Ulam}
\label{_alge_geome_ulam_Subsection_}

The notion of Ulam quasimorphism is not the only 
notion of a quasimorphism considered in the literature.
In \cite{_Hartnick_Schweitzer_}, Hartnick and Schweizer define
an alternative notion of ``a quasimorphism''  $\phi:\; G \arrow H$
as a map such that for any quasimorphism $q:\; H \arrow \R$
its composition with $\phi$ is a quasimorphism
$q\circ \phi:\; G \arrow \R$. When $H$ is 
a hyperbolic group, it has many quasimorphisms
to $\R$ (\cite{_Epstein_Fujiwara_}),
and this condition is quite restrictive. However,
when $H$ is a lattice of high rank, such as $SL(n, \Z)$, 
$n > 2$, (\cite{_Burger_Monod:lattices_,_Burger_Monod:rigidity_}),
all quasimorphisms to $\R$ are equivalent to homomorphisms.
However, the lattices of finite rank satisfy property T
(\cite{_BHV:Kazhdan_T_}), hence they have finite abelinization. 
Therefore, all quasimorphisms on $H$ are bounded; 
for such $H$ the Hartnick-Schweizer condition is quite weak.

\hfill

In \cite{_Fujiwara_Kapovich_}, Fujiwara and Kapovich 
considered other versions of the notion of quasimorphism,
taking values in a non-commutative group $H$; when $H$
is commutative, all these notions are equivalent to the notion of
Ulam quasimorphism.
They define the {\bf algebraic quasimorphisms}
as maps $q:\; G \arrow H$ such that
there exists a compact set $K$ and 
$q(xy) \in K q(x) K q(y) K$ for all $x, y\in G$,
and {\bf geometric quasimorphisms} as ones such that
$q(xy) \in K  q(x) K q(y)$. 

For comparison, Ulam quasimorphisms are maps $q:\; G \arrow H$
which satisfy $q(xy) \in K q(x) q(y)$; this notion is stronger
than the notion of a geometric quasimorphism.

In \cite[Definition 2.2]{_Heuer:bounded_}, N. Heuer introduced a version
of this definition. Given a map $q:\; G \arrow H$,
he defines {\bf the defect} of $q$ as a subgroup generated
by $q(xy)q(y)^{-1} q(y)^{-1}$. In \cite[Proposition 2.3]{_Heuer:bounded_},
Heuer proves that $q$
is an Ulam quasimorphism if and only if the defect is finite.

Since $q(x^{-1}) = q(x)^{-1}$, the condition
$q(xy) \in K q(x) q(y)$ is equivalent to 
$q(xy) \in q(x) q(y) K$. Indeed, 
$q((xy)^{-1}) \in K  q(y^{-1}) q(x^{-1})$
implies $q(xy) \in (K  q(y^{-1}) q(x^{-1}))^{-1} =q(x) q(y) K^{-1}$.
The same argument works for the geometric quasimorphisms
if we also assume that $q(x^{-1}) = q(x)^{-1}$.

Earlier, we defined an equivalence of quasimorphisms:
$q\sim q'$ if there exists a compact $K\subset H$
such that $q(x) \in K q'(x)$ for all $x\in G$.
Interestingly enough, a map which is equivalent 
to an Ulam quasimorphism is no longer an Ulam quasimorphism
(see \ref{_geometric_not_Ulam_Remark_} below).
However, a  map $q':\; G \arrow H$ which is equivalent 
to a geometric quasimorphism $q:\; G \arrow H$  is 
again a geometric quasimorphism:
\begin{align*}
q(xy) &\in Kq(x) Kq(y) \Rightarrow \\
q'(xy) &\in K_1 q(xy) \subset K_1 K q(x) Kq(y)\subset  K_1 KK^{-1}_1 q'(x) KK^{-1}_1 q'(y).
\end{align*}
Here the ``equivalence'' is understood as $q'(x) \in K_1 q(x)$,
and the ``geometric quasimorphism'' condition as
$q(xy) \in Kq(x) Kq(y)$. 

The reason why the ``geometric quasimorphism'' condition
is sometimes more appropriate stems from the following observation.
Let $q:\; G \arrow H$ be an Ulam quasimorphism, and
$d:\; H \times H \arrow \R^{\geq 0}$ a 
left-invariant metric such that all closed
balls in $H$ are compact. Then 
the condition \[ d(q(xy), q(x) q(y)) \leq C\]
is equivalent to being an Ulam quasimorphism.
However, a map $q':\; G \arrow H$ such that 
$d(q(x), q'(x)) \leq C_1$ is no longer an
Ulam quasimorphism, but only a geometric one.
This leads to the following observation.

\hfill

\remark\label{_geometric_not_Ulam_Remark_}
Let $q:\; \Gamma \arrow G$ be a non-constructible
quasimorphism (\ref{_non_constr_Theorem_}) taking values
in a nilpotent Lie group $G$ containing a
cocompact lattice $\Lambda$.
A notion of constructible quasimorphisms
is defined only for Ulam quasimorphisms,
but we can generalize it to geometric quasimorphisms
as follows: a geometric quasimorphism is 
constructible if it is equivalent, in 
the sense of \ref{_Ulam_qm_Definition_}, 
to a constructible Ulam quasimorphism.
Pick a right-invariant metric $d$
on $G$, and let $R$ be the diameter of the 
fundamental domain of $\Lambda$ acting on $G$.
For each $q(x)\in G$, let us choose a closest element 
$q'(x)\in\Lambda$. Then $d(q(x), q'(x))< R$, hence
the map $q':\; \Gamma \arrow \Lambda$ is equivalent
to $q$, and $q'$ is a non-constructible 
geometric quasimorphism.  This gives an example of a geometric quasimorphism
which is non-constructible, even when the target
group is discrete. By \cite{_Fujiwara_Kapovich_},
any Ulam quasimorphism $q'':\; \Gamma \arrow \Lambda$
into a discrete group is constructible, hence $q'$
is not equivalent to any Ulam quasimorphism.




%
%
%
%
%
%

\section{HBG-quasimorphisms with prescribed values}

\subsection{Connections with prescribed holonomy}

We prove the following preliminary lemma, which will be
used later in this section.

\hfill

\lemma\label{_prescribed_conne_Lemma_}
Let $G$ be a connected Lie group, $\g$ its Lie algebra,
and $P$ a trivial $G$-bundle on
an interval $[0, 1]$. Fix an element $g\in G$.
Denote by $\nabla_0$ the trivial connection on $P$.
Then there exists a $\g$-valued 1-form $A$ with
compact support, such that the holonomy $\Hol(\nabla)$
of the connection $\nabla:=\nabla_0 + A$ is equal to $g$.

\hfill

\proof Write $A$ as $a(t)dt$, where $a\in \g$
and $dt$ is the standard 1-form on $[0,1]$.
Then $\Hol(\nabla) = \int_0^1 a(t)dt$.
Since $G$ is connected, we can connect
the unit $1_G$ to $g$ by a path $\gamma:\; [0,1] \arrow G$.
Reparametrizing $\gamma$, we may assume
that $\gamma$ is constant in a small neighbourhood of $0$
and of $1$.
By Newton-Leibnitz formula, 
$\int_0^1 (\gamma(t)^{-1})^*\dot \gamma dt=g$.
Setting $a(t) := (\gamma(t)^{-1})^*\dot \gamma$,
we obtain a connection form which satisfies
$\Hol(\nabla)=\int_0^1 (\gamma(t)^{-1})^*\dot \gamma dt=g$.
Since $\gamma$ is constant in a neighbourhood of $0$ and
of $1$,
the form $a(t) dt$ has compact support.
\endproof

\subsection{Constructing the HBG-quasimorphisms}
\label{_constructing_Subsection_}

In this section, we prove \ref{_existence_Ulam_BG_Theorem_}.
We repeat its statement for convenience.

\hfill

\noindent
{\bf \ref{_existence_Ulam_BG_Theorem_}:\ }
Let $M$ be a closed manifold of strictly negative
curvature, $G$ a non-abelian connected Lie group, and 
$x_1, ..., x_n\in \Gamma:= \pi_1(M)$
a collection of primitive elements satisfying $x_i \not\sim_\vce x_j$
for all $i\neq j$.\footnote{Since $\dim M>1$, there
exist infinitely many such elements,
as follows from \ref{_conj_classes_infinitely_many_Corollary_}.}
Take any collection of elements $g_i\in G$, with $i= 1, ..., n$.
Then 
\begin{description}
\item[(i)]
There exists a connection $\nabla$ on a trivial 
principal bundle $P$ such that the corresponding 
HBG-quasimorphism $q_\nabla:\; \Gamma \arrow G$
takes $x_i$ to $g_i$, $i=1, 2, ..., n$.
\item[(ii)]
Moreover, for any countable family of 
elements $z_j \in \Gamma$ satisfying $z_j \not\sim_\vce x_i$
for all $i, j$, the connection 
$\nabla$ can be chosen in such a way that 
$q_\nabla(z_j)$ is not central in $G$ for all $z_j$.
\end{description}

\pstep
For each primitive element $x\in \pi_1(M)$, denote
by $F_x$ the minimal free geodesic loop representing $x$ 
(\ref{_free_geodesic_unique_Proposition_}).
This geodesic loop is unique and determines $x$ up to conjugation.
Fix a point $p$ on $M$, and a point $p_i$ on each $F_{x_i}$.
Let $\gamma_{x_i}$ be the piecewise smooth loop connecting 
$p$ to $p_i$ by a geodesic segment $\nu_{p_i, F_{x_i}}$, going around $F_{x_i}$, and back
from $p_i$ to $p$ by the same geodesic segment $\nu_{p_i, F_{x_i}}$ reversed.
For each connection $\nabla$ on $P$, the
corresponding HBG-quasimorphism
takes $x_i$ to the holonomy of $\nabla$ along $\gamma_{x_i}$.

We are going to choose a connection $\nabla$ on $P$
such that the holonomy of $\nabla$ along $\gamma_{x_i}$
is equal to $g_i$.

\hfill

{\bf Step 2:}
Since $x_i \not \sim_\vce x_j$, 
none of the loops $F_{x_1}, ..., F_{x_n}$ is contained in another
of these loops. This is actually the only reason why we care
about the VCE equivalence. Fix an open set $B_{x_i}$ containing a segment
of $F_{x_i}$ and not intersecting the rest of the loops.
We can choose $B_{x_i}$ in such a way that
it does not intersect the geodesic segment
connecting $p$ to $p_i$ (Step 1).
Denote by $\nabla_0$ the trivial connection
on $P$. Using \ref{_prescribed_conne_Lemma_}, 
we can construct a connection 1-form 
on each open set $B_{x_i}$ in such a way that
the holonomy of the corresponding connection
 along $\gamma_i\cap B_{x_i}$
is equal to any given element of $G$. 

This gives a new connection on $P$, equal
to the old one outside of $B_{x_i}$, with prescribed
holonomy $\goth H$ on $F_{x_i}$. The holonomy of this connection 
along $\nu_{p_i,F_{x_i}} \circ F_{x_i}  \circ \nu_{p_i,F_{x_i}}^{-1}$
is $\goth S \goth H \goth S^{-1}$, where
$\goth S:\; P\restrict {p}\arrow P\restrict{p_i}$
is the holonomy of this connection along $\nu_{p_i,F_{x_i}}$.
Since the map $\goth H \mapsto \goth S \goth H \goth S^{-1}$
is bijective, by an appropriate choice of 
the connection on $B_{x_i}$, we can obtain any element of $G$
as holonomy of the connection along $\nu_{p_i,F_{x_i}} \circ F_{x_i}  \circ \nu_{p_i,F_{x_i}}^{-1}$.

Using a partition of unity,
we can glue this 1-form to the connection 
form in $\nabla_0$, obtaining another connection
which is equal to $\nabla_0$ outside of $B_{x_i}$,
and has prescribed holonomy on $\gamma_i\cap B_{x_i}$.
This allows us to
modify $\nabla_0$ on each open set $B_{x_i}$ in such a way that
the holonomy of $\nabla$ along $\gamma_i\cap B_{x_i}$
is equal to $g_i$. 

We built a connection $\nabla$ which satisfies
\ref{_existence_Ulam_BG_Theorem_} (i).
To finish the proof of \ref{_existence_Ulam_BG_Theorem_},
it remains to show that $\nabla$ can be chosen 
in such a way that (ii) is satisfied.

\hfill

{\bf Step 3:}
Suppose that $\nabla$ is chosen in such a way that it
satisfies \ref{_existence_Ulam_BG_Theorem_} (i).
We modify the connection $\nabla$ on a sequence of small open sets 
$D_1, ..., D_k$ with each $D_l$ intersecting $\gamma_{z_l}$ and not
intersecting $\gamma_{x_1}, \gamma_{x_2}, ..., \gamma_{x_n}$. The result of these 
successive modifications is a connection denoted $\nabla_k$.

We choose $\nabla_k$ in such a way that 
$q_{\nabla_k}(z_i)$ is not central
for $i = 1, ..., k$. The passage from $\nabla_{k-1}$ to $\nabla_k$
is expressed by $\nabla_k = \nabla_{k-1} + \theta_k$,
where $\theta_k$ is a $\g$-valued 1-form with support in $D_k$. 

We chose the 1-form $\theta_k$
very small, in such a way that the series
$\nabla + \sum \theta_i$ converges to a connection $\tilde \nabla$.
By construction, the holonomy of $\tilde \nabla$ along
$\gamma_{z_k}$ is equal to the holonomy of $\nabla_k$,
hence $q_{\tilde \nabla}(z_k)$ is not central
whenever $q_{\nabla_k}(z_k)=q_{\tilde \nabla}(z_k)$  is not central.

The $\g$-valued 1-form $\theta_k$ is chosen, on each step, 
using the same argument as in the proof of \ref{_prescribed_conne_Lemma_}.
Using this lemma, the 1-form $\theta_k$ is determined
by the value of $q_{\nabla_k}(z_k)$, which can
be chosen arbitrarily small, provided that
$q_{\nabla_k}(z_k)$ is not central. 
Then the series $\nabla + \sum_{i=1}^\infty \theta_i$
converges to a connection $\tilde \nabla$ which satisfies 
\ref{_existence_Ulam_BG_Theorem_} (ii).
Since the support of the form $\sum _{i=1}^\infty\theta_i$
does not intersect the paths $\gamma_{x_1}, ..., \gamma_{x_n}$,
we have $q_{\nabla}(x_i)= q_{\tilde \nabla}(x_i)=g_i$,
hence $\tilde \nabla$ satisfies  both \ref{_existence_Ulam_BG_Theorem_} (i)
and \ref{_existence_Ulam_BG_Theorem_} (ii).
\endproof

\subsection{$\epsilon$-representations}
\label{_epsi_reps_Subsection_}

As an application of the construction given in
Subsection \ref{_constructing_Subsection_}, we 
generalize Kazhdan's theorem \cite[Theorem 2]{_Kazhdan:epsilon_}
to the fundamental group of an arbitrary closed manifold
of strictly negative sectional curvature. 
We define
the notion of an $\epsilon$-representation
and a $\delta$-approximated $\epsilon$-representation
as follows (compare with \cite{_Kazhdan:epsilon_}
and Subsection \ref{_approximation_Subsection_}).

Let $G$ be a topological group.
Recall that {\bf Chabauty topology} on the set
of closed subgroups $C(G)$ is defined by the base of neighbourhoods,
$W_U(\Gamma)$ given for each open subset $U\subset G$
and a subgroup $\Gamma\subset G$
\[
W_U(\Gamma)=\{ \Gamma' \in C(G)\ \ |\ \ U \cdot \Gamma' \supset \Gamma \text{\ \  and\ \ } U \cdot \Gamma \supset \Gamma'\}.
\]
If the topology of $G$ is induced by a left-invariant
metric, the Chabauty topology is induced by 
the Hausdorff metric on the space of subsets of $G$
(\ref{_Hausdorff_distance_Definition_}).

By \cite[p. 474]{_Fell:Hausdorff_}, or 
\cite[Proposition 2 (vi)]{_delaHarpe:Chabauty_},
the space $C(G)$ is compact. 

Recall that a group $G$ {\bf has no small subgroups}
(\cite[Exercise 1.5.6, Corollary 1.5.8]{_Tao:Hilbert_5_})
if there exists a neighbourhood of identity
which does not contain a non-trivial subgroup;
this property is clearly satisfied by any Lie group.
When $G$ has no small subgroups, the space of
$C_0(G)$ closed non-trivial subgroups is also compact.

Clearly, the diameter of a subset of a metric space $M$
is a continuous function in the topology on $2^M$
induced by the Hausdorff metric. Since $C_0(G)$ is compact,
the diameter of non-trivial subgroups in $G$
is bounded from below for any metric on $G$.

Fix a left-invariant
metric on $G$ such that any non-trivial subgroup
has diameter at least $1/3$. 
We motivate the choice of the constant $1/3$ as follows.
For compact group $G$, we can always choose a bi-invariant metric;
its geodesics are translations of one-parametric
subgroups. We normalize the metric such that
the diameter of any compact one-parametric
subgroup is bounded from below by $1/2$.

For a bi-invariant Riemannian metric on a Lie group,
the geodesics are obtained by translation of 
one-parametric group: the Lie-theoretic exponenial
map coincides with the Riemannian. Every compact
Lie group admits a bi-invariant Riemannian metric.

Since every two points are connected
by a geodesic, every element of a compact group 
is an exponent of an element of the Lie algebra,
and every finite cyclic subgroup belongs to a
circle subgroup. If we normalize the metric such that
the diameter of each circle subgroup 
is bounded from below by $1/2$, then the
diameter of each finite subgroup is 
at least $1/3$, which is realized by
$\Z/3\Z$. In particular, of we choose
the bi-invariant metric on the group $SU(2)=S^3$
such that $\diam (SU(2)) = 1/2$ (this is 
equivalent to each meridian circle being of 
length 1), the bound $1/3$ is realized.

{\bf An $\epsilon$-representation
of a group $\Gamma$} is a map $q:\; \Gamma \arrow G$
such that $d(q(x)q(y), q(xy)) < \epsilon$ for any
$x, y \in \Gamma$. An $\epsilon$-representation
can be {\bf $\delta$-approximated by a representation}
if there exists a representation
 $\rho:\; \Gamma \arrow G$ such that
$d (\rho(x), q(x))< \delta$ for all $x\in \Gamma$.

\hfill

\theorem\label{_no_approx_epsilon_rep_Theorem_}
Let $M$ be a closed manifold
of strictly negative sectional curvature, $G$ a positive-dimensional
connected Lie group, and 
$P$ a trivial principal $G$-bundle. For any connection $\nabla$
in $P$, let $q_\nabla:\; \pi_1(M) \arrow G$ denote
the corresponding HBG-quasimorphism
(Subsection \ref{_homo_BG_Subsection_}). Choose a
left-invariant metric on $G$ such that
the diameter of any closed subgroup is at least $1/3$.
Then for each $\epsilon>0$,
there exists a connection $\nabla$ such that
$q_\nabla$ is an $\epsilon$-representation
which cannot be $1/3$-approximated by 
a representation.

\hfill

\pstep 
Let $a_1, ..., a_n$ be the generators of $\pi_1(M)$.
Using \ref{_conj_classes_infinitely_many_Corollary_}, we find
$b \in \pi_1(M)$ such that $b \not\sim_\vce a_i$.
Choose an open set $U_{b}$ which intersects
$F_b$ and does not intersect $F_{a_1}, F_{a_2}, ...$,   
$\nu(p, p_i)$, $i=1, ..., n$, where $\nu(p, p_i)$ 
 are geodesics connecting the
marked point $p$ with $F_{a_i}$.
Let $p_b$ be the chosen fixed point on $F_{b}$,
and $\nabla_0$ be the trivial connection on a 
trivial $G$-bundle $P$.
Choose a connection $\nabla_0+\theta$ which is
trivial outside of $U_{b}$. 
Consider the path $\gamma=\nu(p, p_b) F_{b} \nu(p, p_b)^{-1}$
as a map from $[0,1]$ to $M$. Then 
the holonomy of $\nabla+ \theta$ 
along $\gamma$ is equal to $\int_0^1 \gamma^*(\theta)$.
For any non-trivial $g= e^{u}\in G$ in the image of the exponential 
map $\Lie(G)\to G$, 
we can choose a 1-form $\theta$ on $M$ such that 
$\gamma^*(\theta) = f(t) u dt$, for some function $f:\; [0, 1] \arrow \R$
with compact support. Choosing $f$ such that
$g= \int_0^1 f(t) u dt$, we obtain a 1-form $\theta$
with values in the Lie algebra of $G$ such that
the holonomy of $\nabla_0 + \theta$ along $\gamma$
is equal to $g$. Moreover, replacing $\theta$ with $\frac 1 m \theta$,
we obtain a connection with holonomy $e^{\frac 1 m u} = g^{1/m}$.

For $m$ sufficiently large, this would give
an $\epsilon$-representation $q_\nabla$
such that $q_{\nabla}(a_i) =1_G$,
and $q_\nabla(b)^m =g$.

\hfill

{\bf Step 2:} It remains to show that
the $\epsilon$-representation $q_\nabla$ 
cannot be $1/3$-aproximated by a representation $\rho$.
By contradiction, assume that $q_\nabla$ is 
$1/3$-aproximated by a representation $\rho:\; \pi_1(M) \arrow G$.
Since $d(\rho(a_i^n), q_\nabla(a_i)^n) < 1/3$,
the closure of a subgroup of $G$ generated by $\rho(a_i)$
has diameter less than $1/3$. Since the diameter
of non-trivial subgroups of $G$ is $\geq 1/3$,
this implies that $\rho(a_i)=1_G$,
and $\rho$ is trivial. Therefore, $\rho(b)=1_G$.
However, for all $n\in\Z$, we have
$d(\rho(b)^n, q_\nabla(b)^n) < 1/3$, because $q_\nabla$
is $1/3$-approximated by $\rho$.
Then the diameter of the subgroup
of $G$ generated by $g=q_\nabla(b)$ is
less than $1/3$, which again implies that
$g=1_G$, leading to contradiction.
\endproof

\hfill

{\bf Achnowledgements:} We are very grateful to anonymous referees,
who provided an excellent commentary to our work and gave many 
insightful suggestions and ideas. We are also grateful to IMPA
and Ben Gurion University of the Negev who provided us with
an excellent research atmosphere and working conditions.

\hfill

{\scriptsize

\noindent {\sc Michael Brandenbursky\\
{\sc Ben Gurion University of the Negev\\ 
		Beer Sheva, Israel} \\
{	\tt brandens@bgu.ac.il}}

\hfill

\noindent
{\sc Misha Verbitsky\\
{\sc Instituto Nacional de Matem\'atica Pura e
	Aplicada (IMPA) \\ Estrada Dona Castorina, 110\\
	Jardim Bot\^anico, CEP 22460-320\\
	Rio de Janeiro, RJ - Brasil }\\
also:\\
Laboratory of Algebraic Geometry, \\
Faculty of Mathematics, National Research University HSE,\\
6 Usacheva Str. Moscow, Russia}\\
{\tt verbit@impa.br}}


{\small
\begin{thebibliography}{GMP}


\bibitem[BV1]{_BV:quasimorphisms_orig_}
Michael Brandenbursky, Misha Verbitsky,
{\em Non-commutative Barge-Ghys quasimorphisms},
Int. Math. Res. Not. IMRN 2024, no. 15, 11135-11158.


\bibitem[BV2]{_BV:quasimorphisms_new_}
Michael Brandenbursky, Misha Verbitsky,
{\em Barge-Ghys quasimorphisms
taking values in semi-direct products}, 
\url{http://verbit.ru/TeX/New-quasimorphisms.pdf} and \url{https://www.math.bgu.ac.il/~brandens/New-quasimorphisms.pdf}

\bibitem[DFHM]{_DFHM:quasimorphisms_to_linear_alg_grps_}
Sami Douba, Francesco Fournier-Facio, Sam Hughes, Simon Machado
{\em Quasihomomorphisms to real algebraic groups},
arXiv:2601.22411



\bibitem[Kaz]{_Kazhdan:epsilon_}
D. Kazhdan,
{\em On $\epsilon$-representations},
Israel Journal of Mathematics, vol. {\bf 43}, pages 315-323 (1982),

\bibitem[RW]{_RW:curvature_holonomy_} 
 H. Reckziegel, E. Wilhelmus, 
{\em How the Curvature Generates the Holonomy of a
  Connection in an Arbitrary Fibre Bundle},
Results in Mathematics volume 49, pages 339-359 (2006).



\bibitem[Y]{_Deane_Young:holonomy_}
D. Yang,
{\em Holonomy equals curvature}, 
\url{https://cims.nyu.edu/~yangd/papers/holonomy.pdf}


\end{thebibliography}
}


{\small
\begin{thebibliography}{GMP}

\bibitem[BG]{_Barge_Ghys:bounded_}
J. Barge, \'E. Ghys, 
{\em Surfaces et cohomologie born\'ee,} 
 Invent. Math. 92 (1988), no. 3, 509-526.

\bibitem[BB]{_Bartholdi_Bogopolski_}
L. Bartholdi and O. Bogopolski, {\em On abstract commensurators of groups},
 J. Group Theory, 13(6):903-922, 2010.

\bibitem[BFMSS]{_BFMSS:bounded_}
L. Battista, S. Francaviglia, M. Moraschini, F. Sarti, and
A. Savini, {\em Bounded cohomology classes of
exact forms}, Proc. Amer. Math. Soc., To
appear. arXiv:2211.16125, 2022.



\bibitem[B]{_Bavard:scl_}
C. Bavard, 
{\em Longueur stable des commutateurs},
Enseign. Math. (2) 37 (1991), no. 1-2, 109-150.

\bibitem[BHV]{_BHV:Kazhdan_T_}
B. Bekka, P. de la Harpe, A. Valette,
{\em Kazhdan's Property (T)},
Volume 11 of New Mathematical Monographs,
Cambridge University Press, 2008.


\bibitem[BF]{_Bestvina_Fujiwara_}
M. Bestvina and K. Fujiwara
{\em Bounded cohomology of subgroups of mapping class groups}, 
Geom. Topol., 6, pp 69-89, 2002.

\bibitem[BM]{_Brandenbursky_Marcinkowski_}
M. Brandenbursky, M. Marcinkowski,
{\em Entropy and quasimorphisms,}
Source: Journal of Modern Dynamics (2019), Vol. 15, p 143-163.


\bibitem[BH]{_Bridson_Haefliger:book_}
M. R. Bridson, A. Haefliger, 
{\em Metric spaces of non-positive curvature},
Grundlehren der mathematischen Wissenschaften, 
319. Springer-Verlag, Berlin, 1999. 

\bibitem[Br]{_Brooks:bounded_}
R. Brooks,
{\em Some remarks on bounded cohomology},
Riemann surfaces and related topics: Proceedings of the 1978 Stony Brook Conference (State Univ. New York, Stony Brook, N.Y., 1978), pp. 53-63,
Ann. of Math. Stud., 97, Princeton Univ. Press, Princeton, N.J., 1981.

\bibitem[BBI]{_BBI_}
D. Burago, Yu. Burago,  S. Ivanov,  
{\em A course in metric geometry},
Graduate Studies in Mathematics, 33. American 
Mathematical Society, Providence, RI, 2001.


\bibitem[BM1]{_Burger_Monod:lattices_}
M. Burger and N. Monod, 
{\em Bounded cohomology of lattices in higher rank Lie groups},
J. Eur. Math. Soc. (JEMS) 1 (2) (1999), 199-235.

\bibitem[BM2]{_Burger_Monod:rigidity_}
 Burger, M., and Monod, N., 
{\em Continuous bounded cohomology and applications to rigidity
theory}, Geom. Funct. Anal. 12 (2) (2002), 219-280.



\bibitem[BOT]{_BOT:Ulam_stability_}
M. Burger, N. Ozawa, A. Thom, 
{\em On Ulam stability,}
Israel J. Math. 193 (2013), no. 1, 109-129. 

\bibitem[C]{_Calegari:scl_}
D. Calegari,  {\em scl,}
MSJ Memoirs, 20. Mathematical Society of Japan, Tokyo, 2009. 


\bibitem[CG]{_Corwin_Greenleaf_}
L. Corwin, F. P. Greenleaf, {\em Representations of Nilpotent Lie Groups
and Their Applications. Part I, Basic Theory and Examples,} 
Cambridge Univ. Press, Cambridge, UK, 1990

\bibitem[DGO]{_Dahmani_Guirardel_Osin_}
F. Dahmani, V. Guirardel, and D. Osin. 
{\em Hyperbolically embedded subgroups and rotating families in groups acting on hyperbolic spaces.}
Mem. Amer. Math. Soc., 245(1156):v+152, 2017.

\bibitem[dH]{_delaHarpe:Chabauty_}
P. de la Harpe, {\em 
Spaces of closed subgroups of locally compact groups}, 
2008, arXiv:0807.2030.


\bibitem[EF]{_Epstein_Fujiwara_}
D. B. A. Epstein, K. Fujiwara, 
{\em The second bounded cohomology of word hyperbolic groups,} 
Topology 36, (1997), 1275-1289. 


\bibitem[EP]{_EP:Calabi_qm_}
M. Entov,  L. Polterovich, 
{\em Calabi quasimorphism and quantum homology},
Int. Math. Res. Not. 2003, no. 30, 1635-1676. 

\bibitem[F]{_Fell:Hausdorff_}
J. M. G. Fell, {\em 
A Hausdorff topology for the closed subsets of a locally compact non Hausdorff
space,} Proc. Amer. Math. Soc. 13 (1962), 472-476.


\bibitem[FPS]{_Frigerio_Pozzetti_Sisto_}
R. Frigerio, M. B. Pozzetti, and A. Sisto., 
{\em Extending higher-dimensional quasi-cocycles.}
J.Topol., 8(4):1123-1155, 2015. 

\bibitem[Fri]{_Frigerio:bounded_}
R. Frigerio,
{\em Bounded Cohomology of Discrete Groups},
 Mathematical Surveys and Monographs 227, AMS, 2017.

\bibitem[FK]{_Fujiwara_Kapovich_}
K. Fujiwara, M. Kapovich, 
{\em On quasihomomorphisms with noncommutative targets,}
Geom. Funct. Anal. 26 (2016), no. 2, 478-519.

\bibitem[GG]{_Gambaudo_Ghys_}
J.-M. Gambaudo, \'E. Ghys, 
{\em Commutators and diffeomorphisms of surfaces}.
Ergodic Theory Dynam. Systems 24 (2004), no. 5, 1591-1617.

\bibitem[GLMR]{_GLMR:asymptotic_cohomology_}
L. Glebsky, A. Lubotzky, N. Monod, B. Rangarajan,
{\em Asymptotic Cohomology and Uniform Stability for Lattices in Semisimple Groups},
	arXiv:2301.00476.


\bibitem[G1]{_Gromov:Riemannian_} 
M. Gromov, {\em  Metric structures for Riemannian and
non-Riemannian spaces}, Based on the 1981 French
original. With appendices by M. Katz, P. Pansu and
S. Semmes. Translated from the French by Sean Michael
Bates. Progress in Mathematics, 152. Birkh\"auser Boston,
Inc., Boston, MA, 1999. xx+585 pp.



\bibitem[G2]{_Gromov:Bounded_}
M. Gromov, {\em Volume and bounded
cohomology,} Inst. Hautes \'Etudes Sci. Publ. Math. No. 56
(1982), 5-99 (1983).

\bibitem[HO]{_Hull_Osin_}
M. Hull, D. Osin. 
{\em Induced quasicocycles on groups with hyperbolically embedded subgroups}, 
Algebr. Geom. Topol., 13 (5), pp. 2635-2665, 2013.

\bibitem[HS]{_Hartnick_Schweitzer_}
T. Hartnick, P. Schweitzer,
{\em On quasioutomorphism groups of free groups and their
  transitivity properties},
Journal of Algebra, Volume 450 (2016), pp. 242-281,
arXiv:1403.2786.

\bibitem[He]{_Heuer:bounded_}
N. Heuer, {\em Low-dimensional bounded cohomology and extensions of groups,}
Math. Scand., 126(1):5-31, 2020.

\bibitem[Ho]{_Hochschild:algebraic_Lie_}
G. Hochschild, {\em Note on Algebraic Lie Algebras},
Proceedings of the American Mathematical Society, Vol. 29, 
No. 1 (Jun., 1971), pp. 10-16






\bibitem[Kaz]{_Kazhdan:epsilon_}
D. Kazhdan,
{\em On $\epsilon$-representations},
Israel Journal of Mathematics, vol. {\bf 43}, pages 315-323 (1982),

\bibitem[Kl]{_Klingenberg:Riemannian_} 
W. P. A. Klingenberg, 
{\em Riemannian Geometry}, second edition (1995), de Gruyter.

\bibitem[KN]{_Kob_Nomizu_} 
S. Kobayashi and K. Nomizu, {\em Foundations of differential geometry,}
Wiley-Interscience, New York (1969).

\bibitem[Ko]{_Koch:generators_}
    H. Koch 
{\em Generator and relation ranks for finite-dimensional nilpotent lie algebras,}
Algebra and Logic volume 16, pages 246-253 (1977)


\bibitem[Mar]{_Marasco:Cup_}
D. Marasco,
{\em Cup product in bounded cohomology of negatively curved manifolds}, Proc. Amer. Math.
Soc., 151(6):2707-2715, 2023.

\bibitem[MO1]{_MO:holonomy_question_}
Users of MathOverflow, \\
{\tiny \url{https://mathoverflow.net/questions/423909/holonomy-bounded-in-terms-of-area-and-the-curvature}} \\

\bibitem[MO2]{_MO:primitive_question_}
Y. Cornuiler, S. Nead, Mathoverflow question\\
{\tiny \url{https://mathoverflow.net/questions/423909/holonomy-bounded-in-terms-of-area-and-the-curvature}}\\


\bibitem[M]{_Mineyev:bounded_}
I. Mineyev, {\em Straightening and bounded cohomology of hyperbolic groups},
Geom. Funct. Anal.,
11(4):807-839, 2001.

\bibitem[MW]{_Murphy_Wolhelm:no_geode_subma_}
T. Murphy,  F. Wilhelm, 
{\em Random manifolds have no totally geodesic submanifolds}
Michigan Math. J. 68 (2019), no. 2, 323-335. 

\bibitem[vN]{_von_Neumann:Mechanik_}
J. von Neumann, {\em Beweis des Ergodensatzes und des 
H-Theorems in der neuen Mechanik}, Z. Physik, 57:30-70, 1929.

\bibitem[RW]{_RW:curvature_holonomy_} 
 H. Reckziegel, E. Wilhelmus, 
{\em How the Curvature Generates the Holonomy of a
  Connection in an Arbitrary Fibre Bundle},
Results in Mathematics volume 49, pages 339-359 (2006).

\bibitem[OV]{_Ornea_Verbitsky:Book_} 
L. Ornea, M. Verbitsky, {\em Principles
  of Locally Conformally K\"ahler Geometry},
  arXiv:2208.07188.

\bibitem[O]{_Osin_}
D. Osin,
{\em Acylindrically hyperbolic groups.}
Trans. Amer. Math. Soc., 368(2):851-888, 2016.

\bibitem[P]{_Petersen:Riemannian_} 
P. Petersen,
{\em Riemannian Geometry}, Third Edition, Graduate Texts in Mathematics,
Springer Verlag.

\bibitem[PR]{_RP:functional_}
L. Polterovich, D. Rosen, {\em Function theory on
 symplectic manifolds}, CRM Monograph Series, 34, AMS, Providence, RI, 2014.

\bibitem[R]{_Rolli:free_quasimorphisms_}
P. Rolli,
{\em Quasi-morphisms on free groups}, 
arXiv preprint arXiv:0911.4234, 2009.


\bibitem[Sh]{_Shelukhin:thesis_}
E. Shelukhin, {\em The Action homomorphism,
  quasimorphisms and moment maps on the space of
  compatible almost complex structures},
Commentarii Mathematici Helvetici (2014), 
Volume 89, Issue 1, 2014, 69--123.


\bibitem[St]{_Sternberg:lectures_}
S. Sternberg 
{\em Lectures on differential geometry},
(1964) New York: Chelsea (1093).


\bibitem[Tao]{_Tao:Hilbert_5_}
T. Tao, {\em Hilbert's fifth problem and related topics}, 
Graduate Studies in Mathematics, AMS, 2014.

\bibitem[Tu]{_Turing:Lie_}
A. M. Turing, {\em Finite approximations to Lie groups}, Ann. of
Math. (2), 39(1):105-111, 1938.


\bibitem[U]{_Ulam:stability_}
S. M. Ulam, 
{\em A collection of mathematical problems,}
Interscience Tracts in Pure and Applied Mathematics, 
no. 8 Interscience Publishers, New York-London 1960.

\bibitem[W]{_Wang:597_}
Zh. Wang, {\em Chapter 2: Lattices and nilmanifolds,}
lecture notes to ``MATH 597: Dynamical Systems on
Nilmanifolds'', Penn State University, Fall 2018,
{\scriptsize \url{http://www.personal.psu.edu/zxw14/MATH597/Chapter2.pdf}}




\bibitem[Y]{_Deane_Young:holonomy_}
D. Yang,
{\em Holonomy equals curvature}, 
\url{https://cims.nyu.edu/~yangd/papers/holonomy.pdf}


\end{thebibliography}
}
\end{document}